\documentclass[numbook, envcountsect, envcountsame, envcountreset, runningheads, smallextended]{svjour3}
\usepackage{amssymb}
\usepackage{latexsym}
\usepackage{color}
\usepackage{graphicx}
\usepackage{color}
\usepackage[latin1]{inputenc}
\usepackage[T1]{fontenc}
\usepackage[frenchb]{babel}
\usepackage{lmodern}

\spnewtheorem{contract}{Contract}[section]{\bf}{\it}
\spnewtheorem{assumption}{Assumption}[section]{\bf}{\it}

\def\={\;=\;}
\def\.{\;.}

\def\qeds{\varepsilon}
\def\reff#1{{\rm(\ref{#1})}}
\def \i{1\!\mbox{\rm I}}
\def\1{{\bf 1}}
\def \qed{\hbox{ }\hfill{ ${\cal t}$~\hspace{-5.1mm}~${\cal u}$   } }

\def\normeL2#1{\left\|{#1}\right\|_{L^2}}

\input{tcilatex}
\journalname{Finance and Stochastics}

\begin{document}

\journalname{Finance and Stochastics}
\title{A mathematical treatment of bank monitoring incentives\thanks{%
Research partly supported by the Chair \textit{Financial Risks} of the 
\textit{Risk Foundation} sponsored by Soci\'et\'e G\'en\'erale, the Chair 
\textit{Derivatives of the Future} sponsored by the {F\'ed\'eration Bancaire
Fran\c{c}aise}, and the Chair \textit{Finance and Sustainable Development}
sponsored by EDF and Calyon.} }
\author{ Henri \textsc{Pag\`es} \and Dylan \textsc{Possamaï} }

\titlerunning{A mathematical treatment of bank monitoring incentives}
\authorrunning{Henri Pagès, Dylan Possama\"i}

\institute{H. Pag\`es \at
Banque de France, Reasearch Department, 20 rue du Colonel Driant, 75001,
Paris, France \\\email{henri.pages@banque-france.fr}
\and
D. Possamaï \at CMAP, Ecole Polytechnique, Route de Saclay, 91128, Palaiseau, France \\\email{dylan.possamai@polytechnique.edu.}}

\maketitle

\begin{abstract}
In this paper, we take up the analysis of a principal/agent model with moral
hazard introduced in \cite{pages}, with optimal contracting between
competitive investors and an impatient bank monitoring a pool of long-term
loans subject to Markovian contagion. We provide here a comprehensive
mathematical formulation of the model and show using martingale arguments in
the spirit of Sannikov \cite{san} how the maximization problem with implicit
constraints faced by investors can be reduced to a classical stochastic
control problem. The approach has the advantage of avoiding the more general
techniques based on forward-backward stochastic differential equations
described in \cite{cviz} and leads to a simple recursive system of
Hamilton-Jacobi-Bellman equations. We provide a solution to our problem by a
verification argument and give an explicit description of both the value
function and the optimal contract. Finally, we study the limit case where
the bank is no longer impatient.


\vspace{0.5em}

\end{abstract}

\newpage

\section{Introduction}

Following the seminal contributions of DeMarzo and Fishman \cite{demarzo4}, 
\cite{demarzo5} and Sannikov \cite{san}, there has been a renewed interest
in the mathematical treatment of continuous-time moral hazard models and
their applications. In a typical moral hazard situation, a principal (who
takes the initiative of the contract) is imperfectly informed about the
action of an agent (who accepts or rejects the contract). The goal is to
design a contract that maximizes the utility of the principal while that of
the agent is held to a given level.

\vspace{0.5em} In its whole generality, the mathematical treatment of the
problem can be cast as follows. Agency problems stemming from the agent's
hidden action $a$ limit the utility this agent can get from contracting with
the principal. The optimal contract specifies how these limitations should
be strenghtened or slackened over time as a result of the agent's ongoing
performance. We first have to solve the agent's problem for a given contract 
$c$%
\begin{equation*}
V_{A}(c):=\underset{a}{\sup }\ \mathbb{E}\left[ U_{A}(c,a)\right] ,
\end{equation*}%
where $U_{A}$ is the utility function of the agent. If we assume for
simplicity that there exists a unique optimal action $a(c)$ for any $c$, a
point on the set of constrained Pareto optima can be found by solving the
Principal's stochastic control problem 
\begin{equation*}
V_{P}:=\underset{c}{\sup }\left\{ \mathbb{E}\left[ U_{P}(c,a(c))\right]
+\lambda \mathbb{E}\left[ U_{A}(c,a(c))\right] \right\} ,
\end{equation*}%
where $U_{P}$ is the utility function of the principal and $\lambda $ is the
Lagrange multiplier associated to some reservation utility of the agent.

\vspace{0.5em} Because of the almost limitless choices for $c$, it is
generally assumed that the agent does not have complete control over the
outcomes but instead continuously affects their distribution by choosing
specific actions. This actually means that the agent affects the probability
measure $\mathbb{P}^{a}$ under which the above expectations are taken. This
setting, which will be described more rigorously in the following section,
corresponds to a weak formulation of the stochastic control problem.

\vspace{0.5em} As shown in \cite{cviz}, a general theory can be used to
solve these problems, by means of forward-backward stochastic differential
equations. We show here how recursive, martingale representation-based
techniques proposed by Sannikov \cite{san} can be brought to bear on the
issue to yield explicit solutions that are easy to derive. The paper is a
companion to Pag\`{e}s \cite{pages}, who contributes to the optimal design
of securitization in the presence of banks' impaired incentives to monitor.
It provides a coherent mathematical framework for this problem and lays down
the rigorous foundations for the formal derivations sketched in \cite{pages}.

\vspace{0.5em}Our point is to show that the martingale approach to
contracting can be extended to a Markovian setting, which makes it possible
to relax the assumption of conditional independence between default times.
There are now papers on jump processes as opposed to diffusions, but to the
best of our knowledge they only deal with Poisson risk. In the theory of
repeated games with imperfect monitoring, Abreu et al.\ \cite{abreu} use
Poisson signals to vary the frequency with which actions are taken, and show
that the optimal provision of incentives is markedly different should the
signals be interpreted as \textquotedblleft good\textquotedblright\ or\
\textquotedblleft bad\textquotedblright\ news. Sannikov and Skrzypacz \cite%
{san3} find in a continuous-time setting, with both Brownian information and
Poisson jumps associated with bad news, that deviations from cooperative
behavior can only be punished when the discontinuous information is
revealed. This echoes former results in \cite{san2}, according to which
collusion is impossible under imperfect monitoring with Brownian
information, as the risk of triggering a punishment when no deviation occurs
is large. A model close to ours is Biais et al.\ \cite{biais2}, who deal
with large and unfrequent Poisson losses suffered by a firm that invests in
a stationary environment. They show that an optimal way to restore
incentives when performance is poor is to downsize the project at any time
there is a loss.

\vspace{0.5em}Intensity-based models have been widely used in risk
management. Here we focus on a contagion\ model with \textsl{interacting }%
default intensities. Frey and Backhaus \cite{frey} show how such models can
be conveniently embedded into a Markovian framework. Under this approach, a
Markov chain is defined on the set of all default configurations which, when
names are exchangeable, simply boil down to the portfolio default count, and
default intensities are explicit functions of time and the portfolio default
state, as exemplified by \cite{davis}, \cite{jar} and \cite{yu} in the
credit field. Markov chains have recently been shown to simplify the task of
pricing and hedging credit risk, as in Kraft and Steffensen \cite{kraft} or
Laurent et al.\ \cite{laurent}. They are also useful in our context, as they
provide an alternative way of tackling problems of optimal contracting,%
\footnote{%
With a more complex time dependence, such as the self-exciting Hawkes
formulation of \cite{ait}, \cite{azizpour} or \cite{giesecke}, it may not be
possible to construct a Markov chain describing the jump of each portfolio
constituent. The availability of martingale representation results would
then be questionable.} even though the conditions under which explicit
solutions can be derived are not warranted a priori.

\vspace{0.5em} The rest of the paper is organized as follows. In section \ref%
{Model sec}, we recall the model laid out in \cite{pages}, describe the
contracts and give our main assumptions. In Section \ref{sec.hjb}, we
formally derive a candidate optimal contract by solving the HJB equation
associated to the control problem. We then use a standard verification
argument to show that the candidate solution is indeed the optimal contract
and provide a numerical example. The paper concludes with a short section
devoted to a simple special case.

\section{The model\label{Model sec}}

\subsection{Notations and preliminaries}

\vspace{0.5em} We consider a model with universal risk neutrality in which
time is continuous and indexed by $t\in \left[ 0,\infty \right) $. Without
loss of generality, the risk-free interest rate is taken to be $0$. A bank
has a claim to a pool of $I$ unit loans indexed by $i=1$, $\dots $,$~I$
which are ex ante identical. Each loan is a defaultable perpetuity yielding
cash flow $\mu $ per unit time until it defaults. Once a loan defaults it
gives no further payments. The infinite maturity and no recovery assumptions
are made for tractability.

\vspace{0.5em} Denote by 
\begin{equation*}
N_{t}=\sum_{i=1}^{I}1_{\left\{ \tau ^{i}\leq t\right\} },
\end{equation*}%
the sum of individual loan default indicators, where $\tau ^{i}$ is the
default time of loan~$i$. The current size of the pool\ is $I-N_{t}$. Since
all loans are a priori identical, they can be reindexed in any order after
defaults. The action of the bank consists in deciding at each time $t$
whether it monitors any of the outstanding loans. These actions are
summarized by the functions $e_{t}^{i}$ such that for $1\leq i\leq I-N_{t},\text{ }e_{t}^{i}=1$, if loan $i$ is
monitored at time $t$, and $e_{t}^{i}=0$ otherwise.

\vspace{0.5em}
Non-monitoring renders a private benefit $B>0$ per loan and per unit time to
the bank. The opportunity cost of monitoring is thus proportional to the
number of monitored loans.

\vspace{0.5em} The rate at which loan~$i$ defaults is controlled by the
hazard rate $\alpha _{t}^{i}$ specifying its instantaneous probability of
default conditional on history up to time~$t$. Individual hazard rates are
assumed to depend both on the monitoring choice of the bank and on the size
of the pool. Specifically, we choose to model the hazard rate of a
non-defaulted loan~$i$ at time~$t$ as%
\begin{equation}
\alpha _{t}^{i}=\alpha _{I-N_{t}}\left( 1+(1-e_{t}^{i})\varepsilon \right) ,
\label{hazard eq}
\end{equation}%
where the parameters $\left\{ \alpha _{j}\right\} _{1\leq j\leq I}$
represent individual \textquotedblleft baseline\textquotedblright\ risk
under monitoring when the number of loans is~$j$ and $\varepsilon $ is the
proportional impact of shirking on default risk.

\vspace{0.5em} We define the shirking process $k$ by 
\begin{equation*}
k_{t}=\sum_{i=1}^{I-N_{t}}\left( 1-e_{t}^{i}\right) ,
\end{equation*}%
which represents the number of loans that the bank fails to monitor at time $%
t$. Then, according to \textrm{(\ref{hazard eq})}, aggregate default
intensity is given by 
\begin{equation}
\lambda _{t}^{k}=\alpha _{I-N_{t}}\left( I-N_{t}+\varepsilon k_{t}\right) .
\label{aggregate eq}
\end{equation}

\vspace{0.5em} The bank can fund the pool internally at a cost $r\geq 0$. It
can also raise funds from a competitive investor who values income streams
at the prevailing riskless interest rate of zero. We assume that both the
bank and investors observe the history of defaults and liquidations.

\vspace{0.5em}

\subsection{Description of the contracts}

\vspace{0.5em}Contracts are offered on a take-it-or-leave-it basis by
investors to the bank and agreed upon at time~$0$. They determine how cash
flows are shared and how loans are liquidated, conditionally on past
defaults and liquidations. Without loss of generality, they specify that an
investor receives cash flows from the pool and makes transfers to the bank.
We denote by $D=\left\{ D_{t}\right\} _{t\geq 0}$ the c\`{a}dl\`{a}g,
positive and increasing process describing cumulative transfers from the
investors to the bank, such that 
\begin{equation}
\mathbb{E}^{\mathbb{P}}\left[ D_{\tau }\right] <+\infty ,  \label{integral}
\end{equation}%
where $\tau $ is the liquidation time of the pool and where we assume that $%
D_{0}=0$.

\begin{remark}
For certain interpretations, it will be useful to let $D$ have a jump at
time 0 (cf. Remark \ref{rem.d0}).
\end{remark}

\vspace{0.5em} Let then $H_{t}:=1_{\left\{ t\geq \tau \right\} }$ be the
liquidation indicator of the whole pool. The contract specifies the
probability $\theta _{t}$ with which the pool is maintained given default ($%
dN_{t}=1$), so that at each point in time 
\begin{equation*}
dH_{t}=\left\{ 
\begin{tabular}{ll}
$0$ & with probability $\theta _{t},$ \\ 
$dN_{t}$ & with probability $1-\theta _{t}.$%
\end{tabular}%
\ \right.
\end{equation*}%
With our notations,the hazard rates associated with the default and
liquidation processes $N_{t}$ and $H_{t}$ are $\lambda _{t}^{k}$ and $\left(
1-\theta _{t}\right) \lambda _{t}^{k}$, respectively.

\vspace{0.5em} The contract also specifies when liquidation occurs. We
assume that liquidations can only take the form of the stochastic
liquidation of all loans following immediately default. The above properties
translate into 
\begin{equation*}
\mathbb{P}\left(\tau\in\left\{\tau^1,...,\tau^I\right\}\right)=1,\text{ and }%
\mathbb{P}(\tau=\tau^i|\mathcal{F}_{\tau^i},\tau>\tau^{i-1})=1-\theta_{%
\tau^i}.
\end{equation*}

We summarize the above details of the contracts, which are completely
specified by the choice of $(D,\theta )$. Each infinitesimal time interval $%
\left( t,t+dt\right) $ unfolds as follows:

\begin{itemize}
\item $I-N_{t}$ loans are performing at time $t$.

\item The bank chooses to leave $k_{t}\leq I-N_{t}$ loans unmonitored and
monitors the $I-N_{t}-k_{t}$ others, enjoying private benefits $k_{t}B\,dt$.

\item The investor receives $\left( I-N_{t}\right) \mu \,dt$ from the cash
flows generated by the pool and pays $dDt\geq 0$ as fees to the bank.

\item With probability $\lambda _{t}^{k}\,dt$ defined by (\ref{aggregate eq}%
) there is a default ($dN_{t}=1$).

\item Given default the pool is maintained ($dH_{t}=0$) with probability $%
\theta _{t}$ or liquidated ($dH_{t}=1$) with probability $1-\theta _{t}$.
\end{itemize}

\subsection{Economic assumptions}

\vspace{0.5em}In this section we make some assumptions arising from economic
considerations (see \cite{pages} for details). They are in force throughout
the paper.

\begin{assumption}
\label{assump.mu} 
\begin{equation}
\mu \geq \overline{\alpha }_{I}.
\end{equation}
\end{assumption}

\vspace{0.5em}This condition ensures that monitored loans are profitable
viewed as of time 0.

\begin{assumption}
\label{assump.alpha} We have for all $j\leq I$ 
\begin{equation*}
\frac{r}{\overline{\alpha }_{j}}\leq \frac{\mu \varepsilon -B}{B}\frac{%
\varepsilon }{1+\varepsilon },
\end{equation*}
\end{assumption}

\vspace{0.5em}The condition is related to the efficiency of monitoring and
ensures that the benefits for a non-monitoring bank are not so high that
shirking is socially preferable.

\begin{assumption}
\label{decreasing assumption}Individual default risk is non-decreasing with
past default 
\begin{equation}
\alpha _{j}\leq \alpha _{j-1},\qquad \text{for all }j\leq I.
\label{decreasing cond}
\end{equation}
\end{assumption}

\vspace{0.5em}The condition introduces the possibility of correlated
defaults through a contagion effect, as individual loans' intensity of
default may increase with the arrival of new defaults.

\vspace{0.5em} The expected surplus that can be extracted from the pool of
loans is%
\begin{align}
\nonumber S =\mathbb E\left[ \int_{0}^{\infty }\left( I-N_{t}\right) \mu \,dt\right] -I &=\mathbb E\left[ \sum_{j=1}^{I}\int_{0}^{\infty }j\mu 1_{\left\{ I-N_{t}=j\right\}
}\,dt\right] -I\\
&=I\left( \mu /\overline{\alpha }_{I}-1\right),  \label{first-best}
\end{align}%
which is positive under Assumption \ref{assump.mu}. Indeed, investors could
commit $I$ to the pool, pay the bank $D_{0}=S$ at time 0, and instruct it to
choose $k=0$ until the default of the last loan. They would break even from (%
\ref{first-best}). We assume that the bank's monitoring decision is not
observable. This leads to a dynamic moral hazard problem, as the bank may
choose $k_{t}>0$ down the road to reap the private benefits $k_{t}B$. The
contract $\left( D,\theta \right) $ must use observations on defaults to
give the bank incentives to monitor. We assume that both the bank and
investors can fully commit to such a contract.

\section{Optimal contracting\label{sec.hjb}}

Before going on, let us now describe the stochastic basis on
which we are working. We will always place ourselves on a probability space $%
(\Omega ,\mathbb{F},\mathbb{P})$ on which $N$ is a Poisson process with
intensity $\lambda _{t}^{0}$ (which is defined by \textrm{(\ref{aggregate eq}%
)}) and where $\mathbb{P}$ is the reference probability measure. We denote $(%
\mathcal{F}_{t}^{N})_{t\geq 0})$ the completed natural filtration of $N$ and
by $(\mathcal{G}_{t})_{t\geq 0}$ the minimal filtration containing $(%
\mathcal{F}_{t}^{N})_{t\geq 0})$ and that makes the liquidation time of the
pool $\tau $ a $\mathcal{G}$-stopping time. We note that this filtration
satisfies the usual hypotheses, and therefore we will always consider super
or submartingales in their c\`{a}dl\`{a}g version.

\subsection{Incentive compatibility and limited liability}

As recalled in the introduction, in order to make the problem
tractable, we assume that the monitoring choices of the bank affect the
distribution of the size of the pool. To formalize this, recall that, by
definition, the shirking process $k$ is $\mathcal{G}$-predictable and
bounded. Then, by Girsanov Theorem, we can define a probability measure $%
\mathbb{P}^{k}$ equivalent to $\mathbb{P}$ such that 
\begin{equation*}
N_{t}-\int_{0}^{t}\lambda _{t}^{k}ds,
\end{equation*}%
is a $\mathbb{P}^{k}$-martingale.

\vspace{0.5em} More precisely, we have from Br\'{e}maud \cite{bre} (Chapter
VI, Theorem T3) that on $\mathcal{G}_{t}$ 
\begin{equation*}
\frac{d\mathbb{P}^{k}}{d\mathbb{P}}=Z_{t}^{k},
\end{equation*}%
where $Z^{k}$ is the unique solution of the following SDE 
\begin{equation*}
Z_{t}^{k}=1+\int_{0}^{t}Z_{s^{-}}^{k}\left( \frac{\lambda _{s}^{k}}{\lambda
_{s}^{0}}-1\right) \left( dN_{s}-\lambda _{s}^{0}ds\right) ,\text{ }0\leq
t\leq T,\text{ }\mathbb{P}-a.s.
\end{equation*}%
Then, given a contract $\left( D,\theta \right) $ and a shirking process $k$%
, the bank's expected utility at $t=0$ is given by 
\begin{equation}
u_{0}^{k}(D,\theta ):=\mathbb{E}^{\mathbb{P}^{k}}\left[ \int_{0}^{\tau
}e^{-rt}(dD_{t}+Bk_{t}\,dt)\right] ,  \label{bank utility eq}
\end{equation}%
while that of the investor is 
\begin{equation}
v_{0}^{k}(D,\theta ):=\mathbb{E}^{\mathbb{P}^{k}}\left[ \int_{0}^{\tau
}\left( I-N_{t}\right) \mu \,dt-dD_{t}\right] .  \label{investors utility eq}
\end{equation}

\vspace{0.5em} Following Sannikov \cite{san}, we give now the definition of
an incentive-compatible shirking process.

\begin{definition}
A shirking decision $k$ is incentive-compatible with respect to the contract 
$\left( D,\theta \right) $ if it maximizes \textrm{(\ref{bank utility eq})}.
\end{definition}

\vspace{0.5em} Then, the problem faced by the investors is to design a
contract $\left( D,\theta \right) $ and an incentive-compatible advice on $k$
that maximize their expected discounted payoff, subject to a given
reservation utility for the bank 
\begin{eqnarray}  \label{problem}
v_{I}(u):= &&\underset{(D,\theta )}{\sup }\mathbb{E}^{\mathbb{P}^{k}}\left[
\int_{0}^{\tau }\left( I-N_{t}\right) \mu dt-dD_{t}\right]
\label{Principle agent pb} \\
\text{subject to} &&\mathbb{E}^{\mathbb{P}^{k}}\left[ \int_{0}^{\tau
}e^{-rt}(dD_{t}+Bk_{t}\,dt)\right] \geq u  \notag \\
&&k\text{ incentive-compatible with respect to }\left( D,\theta \right) . 
\notag
\end{eqnarray}%
This allows us to define a first set of admissible contracts for a given
monitoring advice $k$ 
\begin{align}\label{admissible}
\nonumber\mathcal{A}^{k}(x):=\{& (D,\theta ),\text{ }\theta\in[0,1] \text{ is a predictable
process},\ D\text{ is positive, c\`{a}dl\`{a}g,} \\
&  \text{non-decreasing and
satisfies \textrm{(\ref{integral})}, $k$ is incentive-compatible}  \notag \\
& \text{with respect to $\left( D,\theta \right) 
$ and $u_{0}^{k}(D,\theta )\geq x$}\}.
\end{align}

Notice that we will put more restrictions on this set at the end of the
section.

\vspace{0.5em} Using martingale arguments, we now elicit an equivalent
condition for the incentive compatibility of $k$. Consider the bank's
expected lifetime utility, conditional on $\mathcal{G}_{t}$ 
\begin{eqnarray}
U_{t}^{k}(D,\theta ):= &&\mathbb{E}^{\mathbb{P}^{k}}\left[ \int_{0}^{\tau
}e^{-rs}\left( dD_{s}+Bk_{s}ds\right) \Bigm|\mathcal{G}_{t}\right]
\label{lifetime utility eq} \\
&=&\int_{0}^{t\wedge \tau }e^{-rs}\left( dD_{s}+Bk_{s}ds\right)
+e^{-rt}u_{t}^{k}(D,\theta ),  \notag  \label{martingale eq}
\end{eqnarray}%
where $u_{t}^{k}$ is the dynamic version of the bank's continuation utility
defined as 
\begin{equation}
u_{t}^{k}(D,\theta ):=1_{\left\{ t<\tau \right\} }\mathbb{E}^{\mathbb{P}^{k}}%
\left[ \int_{t}^{\tau }e^{-r(s-t)}\left( dD_{s}+Bk_{s}ds\right) \Bigm|%
\mathcal{G}_{t}\right] .  \label{continuation utility eq}
\end{equation}

\vspace{0.5em} Since we are working with the completed natural filtration of
a Poisson process, and since $U_{t}^{k}$ is a $\mathcal{G}_{t}$-martingale
under~$\mathbb{P}^{k}$ and in $L^{1}$ because of the integrability
assumptions we made, the martingale representation theorem for point
processes (see \cite{bre}, Chapter III, Theorems $T9$ and $T17$, and Chapter
VI, Theorems $T2$ and $T3$) implies that there are predictable processes $%
h^{1}$ and $h^{2}$ such that the bank's continuation utility $u^{k}$
satisfies the following \textquotedblleft promise-keeping\textquotedblright\
equation until liquidation occurs 
\begin{equation}
du_{t}^{k}+\left( dD_{t}+Bk_{t}dt\right) =ru_{t}^{k}\,dt-h_{t}^{1}\left(
dN_{t}-\lambda _{t}^{k}\,dt\right) -h_{t}^{2}\left( dH_{t}-(1-\theta
_{t})\lambda _{t}^{k}\,dt\right) ,  \label{PK eq}
\end{equation}%
where the dependence of $h^{1}$ and $h^{2}$ on $k$ has been suppressed for
notational convenience. The introduction of these processes provides a
practical way of characterizing contracts for which a given $k$ is
incentive-compatible, as shown in the following proposition, inspired by
Sannikov \cite{san}. They have the interpretation of \textquotedblleft
penalties\textquotedblright\ weighing down the bank's continuation utility,
the first upon default ($dN_{t}=1$), and the second upon liquidation ($%
dH_{t}=1$).

\begin{proposition}
\label{Sannikov prop}Given a contract $\left( D,\theta \right) $ and a
shirking process $k$, the latter is incentive-compatible if and only if for
all $t\in \left[ 0,\tau \right] $ and for all $i=1$, $\cdots$, $I-N_t$, the
following holds almost-surely, 
\begin{equation}
\left( \frac{B}{\varepsilon \alpha _{I-N_{t}}}-h_{t}^{1}-(1-\theta
_{t})h_{t}^{2}\right) (k_{t}-i)\geq 0.  \label{IC eq}
\end{equation}
\end{proposition}

\vspace{0.5em} {\noindent \textbf{Proof. }} Consider an arbitrary strategy $%
\widehat{k}$ specifying the number of unmonitored loans at any point in time
until liquidation. Let $u_{t}^{k}$ denote the continuation utility in (\ref%
{continuation utility eq}) resulting from the decision to forgo monitoring $%
k $ loans at all times.

Define by%
\begin{equation}
\widehat{U}_{t}=\int_{0}^{t\wedge \tau }e^{-rs}\left( dD_{s}+B\widehat{k}%
_{s}ds\right) +e^{-rt}u_{t}^{k}  \label{eq.U}
\end{equation}%
the lifetime utility of the bank viewed as of time$~t$ if it follows the
strategy $\widehat{k}$ before time~$t$, and plans to switch to $k$
afterwards.

\vspace{0.5em} We have for all $t\in \left[ 0,\tau \right] $ 
\begin{align*}
d\widehat{U}_{t}& =e^{-rt}\left( dD_{t}+B\widehat{k}_{t}dt\right)
+e^{-rt}\left( du_{t}^{k}-ru_{t}^{k}\,dt\right) \\
& =e^{-rt}B(\widehat{k}_{t}-k_{t})\,dt-e^{-rt}\left(
h_{t}^{1}(dN_{t}-\lambda _{t}^{k}\,dt)+h_{t}^{2}(dH_{t}-(1-\theta
_{t})\lambda _{t}^{k}\,dt)\right) \\
& =e^{-rt}\left( B-\alpha _{I-N_{t}}\varepsilon (h_{t}^{1}+(1-\theta
_{t})h_{t}^{2})\right) (\widehat{k}_{t}-k_{t})dt \\
& \hspace{0.9em}-e^{-rt}\left( h_{t}^{1}(dN_{t}-\lambda _{t}^{\widehat{k}%
}\,dt)+h_{t}^{2}(dH_{t}-(1-\theta _{t})\lambda _{t}^{\widehat{k}%
}\,dt)\right) ,
\end{align*}%
where we have used the promise-keeping equation~(\ref{PK eq}) for~ $u^{k}$.
Therefore, the first term on the right-hand side 
\begin{equation*}
e^{-rt}\left( B-\alpha _{I-N_{t}}\varepsilon (h_{t}^{1}+(1-\theta
_{t})h_{t}^{2})\right) (\widehat{k}_{t}-k_{t}),
\end{equation*}%
is the drift of $\widehat{U}$ under $\mathbb{P}^{\widehat{k}}$. Note also
that, by definition, $h^{1}$ and $h^{2}$ are integrable and therefore the
martingale part of $\widehat{U}$ is a true $\mathbb{P}^{\widehat{k}}$%
-martingale.

\vspace{0.5em} $\mathrm{(i)}$ Now assume that \textrm{(\ref{IC eq})} does
not hold on a set of positive measure, and choose $\widehat{k}$ such that it
maximizes the quantity 
\begin{equation*}
\left( B-\alpha _{I-N_{t}}\varepsilon (h_{t}^{1}+(1-\theta
_{t})h_{t}^{2})\right) \widehat{k}_{t},
\end{equation*}%
for all $t$. Then, the drift of $\widehat{U}$ under $\mathbb{P}^{\widehat{k}%
} $ is non-negative and strictly positive on a set of positive measure.
Therefore $\widehat{U}$ is a $\mathbb{P}^{\widehat{k}}$-submartingale. This
implies the existence of a time $t^{\ast }>0$ such that 
\begin{equation*}
\mathbb{E}^{\mathbb{P}^{\widehat{k}}}[\widehat{U}_{t^{\ast }}]>\widehat{U}%
_{0}=u_{0}^{k}.
\end{equation*}%
Therefore, if the agent follows this strategy $\widehat{k}$ until the time $%
t^{\ast }$ and then switches to the strategy $k$, his utility is strictly
greater than the utility obtained from following the strategy $k$ all the
time. This contradicts the fact that the strategy $k$ is
incentive-compatible.

\vspace{0.5em} $\mathrm{(ii)}$ With the same notations as above, assume that 
\textrm{(\ref{IC eq})} holds for the strategy $k$. Then this means that $%
\widehat{U}$ is a $\mathbb{P}^{\widehat{k}}$-supermartingale, regardless of
the choice of strategy $\widehat{k}$. Moreover, since $\widehat{U}$ is
positive (because $D$ is non-decreasing), it has a last element (see Problem 
$3.16$ in \cite{karat} for instance). Then, we have by the optional sampling
Theorem

\begin{equation*}
u_{0}^{k}=\widehat{U}_{0}\geq \mathbb{E}^{\mathbb{P}^{\widehat{k}}}\left[ 
\widehat{U}_{\tau }\right] =u_{0}^{\widehat{k}},
\end{equation*}%
where we used \textrm{(\ref{eq.U})} and the fact that $u_{\tau }^{k}=0$ for
the last inequality. This means that the strategy $k$ maximizes the expected
utility of the agent and is therefore incentive-compatible. \hbox{ }\hfill $%
\Box $

\vspace{0.5em}Under the assumption that monitoring is efficient, we now
focus on contracts that actually deter the bank from shirking, i.e.,
contracts with respect to which $k=0$ is incentive-compatible. In that
particular case, the above Proposition can be simplified as follows.

\begin{corollary}
\label{Sannikov prop2}Given a contract $\left( D,\theta \right) $, $k=0$ is
incentive-compatible if and only if 
\begin{equation}
h_{t}^{1}+(1-\theta _{t})h_{t}^{2}\geq \frac{B}{\varepsilon \alpha _{I-N_{t}}%
},\text{ }t\in \lbrack 0,\tau ],\text{ }\mathbb{P}-a.s.  \label{IC cond}
\end{equation}
\end{corollary}

\vspace{0.5em}

\begin{remark}
Corollary~\ref{Sannikov prop2} states that, given that the pool has $i$
loans outstanding, in order to induce the bank to monitor all loans, the
continuation payoff must drop in expectation by at least the quantity 
\begin{equation*}
b_{i}:=\frac{B}{\varepsilon \alpha _{i}},
\end{equation*}%
following default.
\end{remark}

\vspace{0.5em} In order to specify further our admissible strategies, we
have to put some restrictions on $h^{1}$ and $h^{2}$. First, we assume that
the bank has limited liability. This means that the bank's continuation
utility is bounded from below by $b_{I-N_{t}}$ up to liquidation, since
otherwise the incentive-compatible (\ref{IC cond}) would be violated upon
default. In particular, the limited liability constraint must hold after a
default if the pool is maintained in operation ($dH_{t}=0$), when the drop
in utility is $h^{1}$. This implies that%
\begin{equation}
\text{For all $1\leq i\leq I$, }u_{t^{-}}^{0}-h_{t}^{1}\geq b_{i-1},\text{
on $\left\{ N_{t}=I-i\right\} $.}  \label{limited liability eq}
\end{equation}%
For the second condition, we assume that the bank forfeits any rights to
cash flows once the pool is liquidated. The constraint $u_{\tau }^{0}=0$
implies in turn that at all times 
\begin{equation}
u_{t^{-}}^{0}=h_{t}^{1}+h_{t}^{2},  \label{liquidation eq}
\end{equation}%
since the drop in utility is $h^{1}+h^{2}$ in that case.

\vspace{0.5em} The introduction of the processes $h^{1}$ and $h^{2}$ allows
us to greatly simplify the set of admissible contracts by formulating the
incentive compatibility requirement in terms of explicit conditions. Our set
of admissible strategies is therefore

\begin{align}
\widetilde{\mathcal{A}}^{0}(x):=\{& (D,\theta ,h^{1},h^{2}),\text{ }\theta\in[0,1] 
\text{ is predictable},\ D\text{ is positive, c\`{a}dl\`{a}g,}  \notag
\label{admissible2} \\
&  \text{non-decreasing and
satisfies \textrm{(\ref{integral})}, } \text{$h^{1}$ and $h^{2}$ are predictable,}  \notag \\
& \text{integrable and
satisfy $u_{t^{-}}^{0}-h_{t}^{1}\geq b_{I-N_{t}-1}$}, \ \text{$u_{t^{-}}^{0}=h_{t}^{1}+h_{t}^{2}$,} \notag \\
&  \text{ and $x\leq u_{0}^{0}(D,\theta )$.}%
\}.
\end{align}%
Since $j=1$ is a degenerate special case, it is convenient to treat
monitoring with a single loan first before turning to the general case.

\subsection{Single loan: Constant utility}

\vspace{0.5em}The default of a single loan ends the game. Hence, there is no
room for stochastic liquidation and the processes $\theta $ and $h^{2}$ are
left undefined. The incentive compatibility constraint (\ref{IC cond}) takes
the simpler form $h_{t}^{1}\geq b_{1}$. However, from \textrm{(\ref%
{liquidation eq})}, $u_{t}=h_{t}^{1}$ on $\left\{ t<\tau \right\} $, so the
incentive compability constraint can be rewritten as:%
\begin{equation*}
u_{t}\geq b_{1},\ \quad t<\tau \text{, }\mathbb{P}-a.s.
\end{equation*}%
Note that the limited liability constraint is automatically satisfied and
that (\ref{limited liability eq}) can be disregarded as $u_{\tau }=0$ upon
default. Noting $v_{1}(u)$ the highest value that an investor can achieve
for a given bank's continuation value of$~u$, we have the following result.

\begin{proposition}
\label{prop.un} When $j=1$, the value function is given by $v_{1}(u)=b_{1}-u+%
\overline{v}_{1}$ on $u\geq b_{1}$, where%
\begin{equation*}
\overline{v}_{1}:=\frac{\mu -b_{1}(r+\lambda _{1})}{\lambda _{1}}.
\end{equation*}%
Under the optimal contract, starting from reservation utility $u\geq b_{1}$,
the incentive compatibility constraint binds at all times until default. The
bank receives:

\begin{itemize}
\item An initial lump-sum payment $D_{0}=u-b_{1}$ which brings its
continuation utility back to $b_{1}$,

\item A continuous payment $dD_{t}=b_{1}(r+\lambda _{1})\,dt$ until default.
\end{itemize}
\end{proposition}

The proof is relegated to the Appendix.

\subsection{Reduction to a stochastic control problem and HJB equations}

\vspace{0.5em}Let us now turn to the general case $j\geq 2$. Under condition
(%
\ref{IC cond}%
), $k=0$ is incentive-compatible. That being taken care of, solving for the
optimal contract involves maximizing an investor's expected utility and is
therefore a classical stochastic control problem. Let $v_{j}(u)$ denote the
investor's value function, i.e., the maximum expected utility an investor
can achieve given a pool of size $j$ and a reservation utility $u$ for the
bank. Assume for now that the process $D$ is absolutely continuous with
respect to the Lebesgue measure (we will verify later that the property is
satisfied at the optimum), that is to say 
\begin{equation*}
D_{t}=\int_{0}^{t}\delta _{s}ds.
\end{equation*}

We expect the investor's value function to solve the following
system of HJB equations with initial condition $v_{1}(u)$%
\begin{align}
& \underset{(\delta ,\theta ,h^{1},h^{2})\in \mathcal{C}^{j}}{\sup }\left\{
\left( ru+\lambda _{j}\left( h^{1}+(1-\theta )h^{2}\right) -\delta \right)
v_{j}^{\prime }(u)+j\mu -\delta \right.  \notag  \label{HJBB} \\
& \hspace{6.5em}\left. -\theta \lambda _{j}\left(
v_{j}(u)-v_{j-1}(u-h^{1})\right) -(1-\theta )\lambda _{j}v_{j}(u)\right\} =0,%
\text{ }u\geq b_{j},
\end{align}%
where the $\mathcal{C}^{j}$ are our admissible strategies sets defined by 
\begin{equation*}
\mathcal{C}^{j}:=\left\{\delta \geq 0,%
\text{ }\theta \in \lbrack 0,1],\text{ }h^{1}+(1-\theta )h^{2}\geq b_{j},%
\text{ }u-h^{1}\geq b_{j-1},\text{ }u=h^{1}+h^{2}\right\} .
\end{equation*}

\vspace{0.5em}
\begin{remark}
We will see in the next section that our control problem is singular.
Therefore the above HJB equation \textrm{(\ref{HJBB})} is not exactly the
correct one, and we will consider instead a variational inequality.
\end{remark}

\vspace{0.5em} Given the constraints in the definition of $\mathcal{C}^{j}$,
we reparametrize the problem in terms of the variable $z:=\theta (u-h^{1})$.
This leads to the simpler system of HJB equations for $u\geq b_j$

\begin{equation}
\underset{(\delta ,\theta ,z)\in \widetilde{\mathcal{C}}^{j}}{\sup }\left\{
\left( ru+\lambda _{j}\left( u-z\right) -\delta \right) v_{j}^{\prime
}(u)+j\mu -\delta -\lambda _{j}( v_{j}(u)-\theta v_{j-1}( \frac{z}{
\theta })) \right\} =0,  \label{HJB bis eq}
\end{equation}%
where the constraints become 
\begin{equation*}
\widetilde{\mathcal{C}}^{j}:=\left\{ (\delta ,\theta ,z),\text{ }\delta \geq
0,\text{ }\theta \in \left[ 0,1\wedge \frac{u-b_{j}}{b_{j-1}}\right] ,\text{
and }z\in \lbrack b_{j-1}\theta ,u-b_{j}]\right\} .
\end{equation*}%
Our strategy now is to guess a candidate optimal contract by solving the
above system of HJB equations, and to prove that the conjectured contract is
indeed optimal by means of a verification argument.

\subsubsection{Formal derivation of a candidate optimal contract\label%
{formal derivation}}

\vspace{0.5em} \textbf{Step $\mathrm{(i)}$} Optimizing first with respect to 
$\delta $ yields the following variational inequality for $u>b_{j}$ 
\begin{align}
\nonumber& \min \left\{ -\underset{(\theta ,z)\in \widetilde{\mathcal{B}}^{j}}{\sup }%
\left\{ \left( ru+\lambda _{j}\left( u-z\right) \right) v_{j}^{\prime
}(u)+j\mu -\lambda _{j}\left( v_{j}(u)-\theta v_{j-1}\left( \frac{z}{\theta }%
\right) \right) \right\},\right.\\
&\hspace{3.2em}\left.v_{j}^{\prime }(u)+1\right\} =0. 
\label{trucbank}
\end{align}%
where 
\begin{equation*}
\widetilde{\mathcal{B}}^{j}:=\left\{ (\theta ,z),\text{ }\theta \in \left[
0,1\wedge \frac{u-b_{j}}{b_{j-1}}\right] ,\text{ and }z\in \lbrack
b_{j-1}\theta ,u-b_{j}]\right\} .
\end{equation*}

We assume that all the functions $v_{j}$ are concave (a property which needs
to be verified by our candidate). Then the first derivative of $v_{j}$ is
decreasing. Let us also assume that there exists a level $\gamma _{j}>b_{j}$
(a free boundary) such that: 
\begin{equation*}
v_{j}^{\prime }(\gamma _{j})=-1,\text{ }v_{j}^{\prime }(u)>-1,\text{ for $%
u<\gamma _{j}.$}
\end{equation*}%
Then as long as $u<\gamma _{j}$, $v_{j}$ satisfies the first equation in 
\textrm{(\ref{trucbank})}. Therefore, equation \textrm{(\ref{trucbank})}
tells us that the bank cannot receive cash from investors unless its utility
attains the level$~\gamma _{j}$ (since $\delta =0$ is optimal before that).
We also assume (and will verify) that our candidate satisfy for $u\geq \gamma_j$
\begin{equation*}
\underset{(\theta ,z)\in \widetilde{\mathcal{B}}^{j}}{\sup }\left\{ \left(
ru+\lambda _{j}\left( u-z\right) \right) v_{j}^{\prime }(u)+j\mu -\lambda
_{j}\left( v_{j}(u)-\theta v_{j-1}\left( \frac{z}{\theta }\right) \right)
\right\} \leq 0.
\end{equation*}

This means that $v_{j}$ becomes linear above $\gamma _{j}$, and that the
variational inequality \textrm{(\ref{trucbank})} takes the following simpler form. If $u\in(b_j,\gamma_j]$ then 
\begin{align*}
& \underset{(\theta ,z)\in \widetilde{\mathcal{B}}^{j}}{\sup }\left\{
\left( ru+\lambda _{j}\left( u-z\right) \right) v_{j}^{\prime }(u)+j\mu
-\lambda _{j}\left( v_{j}(u)-\theta v_{j-1}\left( \frac{z}{\theta }\right)
\right) \right\} =0,
\end{align*}
and if $u>\gamma_j$
\begin{align*}
& v_{j}^{\prime }(u)+1=0.
\end{align*}

\vspace{0.5em} In order to choose $\gamma _{j}$, it is natural to require
our solution to be maximal in the sense that for each $u>b_{j}$:%
\begin{equation*}
\gamma _{j}\longrightarrow v_{j}(u),
\end{equation*}%
is maximal at the chosen value of $\gamma _{j}$. Of course, it is not clear
whether such a value exists. This heuristic approach can be proven
rigorously, and that the maximality assumption will be clarified.

\begin{property}
Payments are made to the bank only when its continuation utility reaches a
threshold $\gamma _{j}$ satisfying $v^{\prime }(\gamma _{j})=-1$.
\end{property}

\vspace{0.5em}The economic interpretation is as follows. Under the assumed
concavity of the continuation function, its slope $v_{j}^{\prime }(u_{j})$
is strictly above $-1$ as long as total utility $u+v_{j}(u)$\ fails to be
maximized. Put differently, it is less expensive for investors to allow for
an increase in the bank's continuation payoff, which costs them $%
v_{j}^{\prime }(u_{j})$, than paying the bank right away, which costs them $%
-1$. Compensating the bank with a higher continuation utility rather than
cash implies that payments are deferred ($\delta =0$). If no default incurs,
the bank's continuation payoff $u$ keeps increasing as the time set for the
resumption of payments gets closer and closer. It eventually reaches the
optimum level $\gamma _{j}$ (unless a default interrupts the process), at
which point the bank is paid. To this extent, the state variable $u$ can be
interpreted as a measure of performance, with the derivation showing that
the optimal compensation scheme must be based on performance.

\vspace{0.5em}\textbf{Step $\mathrm{(ii)}$} We next turn to the liquidation
decision. One finds as first-order condition with respect to $\theta $: 
\begin{equation}
v_{j-1}\left( \frac{z}{\theta }\right) -\frac{z}{\theta }v_{j-1}^{\prime
}\left( \frac{z}{\theta }\right) \geq 0.  \label{FOC theta eq}
\end{equation}%
Once again, if $v_{j-1}$ is concave, the above inequality (\ref{FOC theta eq}%
) is always verified. This means that the function 
\begin{equation*}
\theta \longrightarrow \theta v_{j-1}\left( \frac{z}{\theta }\right) ,
\end{equation*}%
is non-decreasing, which implies that the optimal $\theta $ corresponds to
its upper bound. There are then two cases:

\begin{itemize}
\item[(i)] $u\in \left[ b_{j},b_{j}+b_{j-1}\right) $ and $\theta =\left(
u-b_{j}\right) /b_{j-1}$

\item[(ii)] $u\in \left[ b_{j}+b_{j-1},\gamma _{j}\right) $ and $\theta =1$.
\end{itemize}

\begin{property}
Stochastic liquidation takes place in the interval $\left[
b_{j},b_{j}+b_{j-1}\right) $, with controls given by%
\begin{equation*}
\left\{ 
\begin{array}{rcl}
\delta _{t} & = & 0 \\ 
\theta _{t} & = & \left( u_{t}-b_{j}\right) /b_{j-1} \\ 
h_{t}^{1} & = & u_{t}-b_{j-1} \\ 
h_{t}^{2} & = & b_{j-1}%
\end{array}%
\right.
\end{equation*}%
There is no liquidation in the probation interval $\left[ b_{j}+b_{j-1},%
\gamma _{j}\right) $.
\end{property}

\vspace{0.5em} In the first interval, the pool is liquidated with strictly
positive probability following default. Since $b_{j-1}\theta =u-b_{j}$, the
only value left for $z$ is $z=u-b_{j}=b_{j-1}\theta $, from which we derive $%
h^{1}=u-b_{j-1}$ and $h^{2}=b_{j-1}$. Thus, if a default occurs in that
interval, either the bank's continuation utility drops to the minimum
threshold $b_{j-1}$ (with probability $\theta $ given by the position of $u$
in that interval) or the pool is liquidated. This actually ensures that the
incentive compatibility condition (\ref{IC cond}) is met, even though under
continuation the drop in continuation utility $u-b_{j-1}$ is below the
minimum $b_{j}$ required for incentive purposes. In contrast, there is no
liquidation in the interval $\left[ b_{j}+b_{j-1},\gamma _{j}\right) $,
which we refer to as \textquotedblleft probation.\textquotedblright\ It will
be verified that $\gamma _{j}\geq b_{j}+b_{j-1}$, implying that the
stochastic liquidation interval has always a width of $b_{j-1}$.

\vspace{0.5em} \textbf{Step $\mathrm{(iii)}$} Finally consider the decision
regarding $z$. We have seen that, if $u\in \left[ b_{j},b_{j}+b_{j-1}\right) 
$, then $z=u-b_{j}$. On $u\in \left[ b_{j}+b_{j-1},\gamma _{j}\right] $, $%
\theta =1$ and $z$ is constrained in the range $\left[ b_{j-1},u-b_{j}\right]
$. We continue our guess of a candidate solution by assuming that 
\begin{equation}
v_{j-1}^{\prime }(u-b_{j})-v_{j}^{\prime }(u)\geq 0,  \label{FOC z eq}
\end{equation}%
a condition which needs to be verified by the resulting candidate. Since $%
v_{j-1}$ is supposed to be concave, we have for all $z\in \left[
b_{j-1},u-b_{j}\right] $ 
\begin{equation}
v_{j-1}^{\prime }(z)-v_{j}^{\prime }(u)\geq 0.  \label{FOC z eq bis}
\end{equation}%
From this, we obtain that the function $z\longrightarrow -zv_{j}^{\prime
}+v_{j-1}(z)$ is non-decreasing, which in turn implies that the supremum
over $z$ is attained at $u-b_{j}$. This implies in turn that, when $u=\gamma
_{j}$, the dividend payment is $\delta _{j}=ru+\lambda _{j}(u-z)=\lambda
_{j}b_{j}+r\gamma _{j}$.

\begin{property}
The incentive
compatibility constraint (\ref{IC cond}) binds on the interval $\left[ b_{j}+b_{j-1},\gamma _{j}\right] $, and controls are given by%
\begin{equation*}
\left\{ 
\begin{array}{rcl}
\delta _{t} & = & 1_{\left\{ u_{t}=\gamma _{j}\right\} }\left( \lambda
_{j}b_{j}+r\gamma _{j}\right) \\ 
\theta _{t} & = & 1 \\ 
h_{t}^{1} & = & b_{j} \\ 
h_{t}^{2} & = & u_{t}-b_{j}.%
\end{array}%
\right.
\end{equation*}
\end{property}

\vspace{0.5em}Note that the value assigned to $h^{2}$ is irrelevant, as
stochastic liquidation is never carried out in this interval. The idea
behind (\ref{FOC z eq bis}) is that, whatever the choice of $z$, investor
value becomes more sensitive to performance as the bank's continuation
utility takes a cut following default from $u$ to $z$. As long as the
difference is positive, the investor is willing to increase $z$ up to its
maximum $u-b_{j}$, i.e., reduce the penalty to the incentive-compatible
level $h^{1}=b_{j}$. Intuitively, it is costly to impose a higher penalty
than necessary, because it would require that the bank be compensated with a
higher utility growth under probation or with higher payments at the
threshold $\gamma _{j}$, which would reduce investor value. Finally, note
that the dividend flow $\delta $ has two components. The first, $\lambda
_{j}b_{j}=jB/\epsilon $, is proportional to size and can be interpreted as a
monitoring (or servicing) fee. The second, $r\gamma _{j}$, is tuned to the
bank's rate of impatience and can be interpreted as a \textquotedblleft
rent-preserving\textquotedblright\ fee. The performance-based compensation
scheme resembles that obtained in actual securitization arrangements, at
least when the sponsor retains an equity tranche, with some important
differences that are streamlined in \cite{pages}.

\vspace{0.5em} Summarizing all the above formal calculations, we can finally
describe the contract $(D,\theta )$.

\begin{contract}
\label{contract}For given size $j\in \left\{ 1,\dots ,I\right\} $, let the
controls in (\ref{HJBB}) be defined as:%
\begin{align}
& \delta ^{j}(u):=1_{\left\{ u=\gamma _{j}\right\} }(\lambda
_{j}b_{j}+r\gamma _{j})  \notag  \label{contrat} \\
& \theta ^{j}(u):=1_{\left\{ b_{j}\leq u<b_{j}+b_{j-1}\right\} }\left(
u-b_{j}\right) /b_{j-1}+1_{\left\{ b_{j}+b_{j-1}\leq u\leq \gamma
_{j}\right\} }  \notag \\
& h^{1,j}(u):=(u-b_{j-1})1_{\left\{ b_{j}\leq u<b_{j}+b_{j-1}\right\}
}+b_{j}1_{\left\{ b_{j}+b_{j-1}\leq u\leq \gamma _{j}\right\} }  \notag \\
& h^{2,j}(u):=u-h^{1,j}(u).
\end{align}%
The corresponding contract can be described as follows:

\begin{itemize}
\item[$\mathrm{{(i)}}$] Given size $j$, the pool remains in operation (i.e.
there is no liquidation) with one less unit at any time there is a default
in the range $\left[ b_{j}+b_{j-1},\gamma _{j}\right] .$

\item[$\mathrm{{(ii)}}$] The flow of dividend paid to the bank given $j$ is $%
\delta _{t}^{j}=\lambda _{j}b_{j}+r\gamma _{j}$ as long as $u_{t}=\gamma
_{j} $ and no default occurs, where $\delta ^{j}$ is the density of $D$ with
respect to the Lebesgue measure. Otherwise $\delta _{t}=0$.

\item[$\mathrm{{(iii)}}$] Liquidation of the pool occurs with probability $%
\theta _{t}^{j}=\left( u_{t}-b_{j}\right) /b_{j-1}$ at any time there is a
default in the range $\left[ b_{j},b_{j}+b_{j-1}\right) $. If the pool is
maintained, the bank's continuation utility is reset to its minimum $b_{j-1}$
consistent with size $j-1$.
\end{itemize}
\end{contract}

\begin{remark}
\label{rem.d0} If the bank's reservation utility at time $0$ is greater than 
$\gamma _{I}$, then the contract should specify that a transfer is
immediately made to the bank so that its utility is brought back to $\gamma
_{I}$. This means that instead of considering transfers $(D_{t})_{t\geq 0}$
which are absolutely continuous with respect to the Lebesgue measure, we
have to add a Dirac mass at $0$. This can be readily shown from the form of
the value function $v_{I}(u)$, therefore we will not treat it. Notice that
the contract \ref{contract} is clearly in $\widetilde{\mathcal{A}}^{0}(x)$.
\end{remark}

\vspace{0.5em}We end up with the following system of ODEs characterizing the
HJB equations on the interval $[b_{j},\gamma _{j}]$
\begin{align*}
& \left( ru+\lambda _{j}b_{j}\right) v_{j}^{\prime }(u)+j\mu -\lambda
_{j}\left( v_{j}(u)-v_{j-1}(u-b_{j})\right) =0,\ u\in \left(
b_{j}+b_{j-1},\gamma _{j}\right] \\
& \left( ru+\lambda _{j}b_{j}\right) v_{j}^{\prime }(u)+j\mu -\lambda
_{j}\left( v_{j}(u)-\scriptstyle\frac{u-b_{j}}{b_{j-1}}v_{j-1}(b_{j-1})\right) =0,\
u\in \left( b_{j},b_{j}+b_{j-1}\right] .
\end{align*}
We can simplify somewhat the formulation by extending the value function $%
v_{j}$ to the interval $[0,b_{j}]$ as
\begin{equation}
v_{j}(u):=\frac{u}{b_{j}}v_{j}(b_{j}),\ u\in \left[ 0,b_{j}\right] ,
\label{extrapolation eq}
\end{equation}%
and to the interval $(\gamma _{j},+\infty )$ as: 
\begin{equation*}
v_{j}(u):=v_{j}(\gamma _{j})-u+\gamma _{j}.
\end{equation*}%
Then the above system of ODEs becomes:%
\begin{align}
& \left( ru+\lambda _{j}b_{j}\right) v_{j}^{\prime }(u)+j\mu -\lambda
_{j}\left( v_{j}(u)-v_{j-1}(u-b_{j})\right) =0,\ u\in \left(
b_{j},\gamma _{j}\right]  \label{HJB.new} \\
& v_{j}^{\prime }(u)=-1,\text{ }u\geq \gamma _{j}.  \notag
\end{align}
We need to verify that the solution obtained from \textrm{(\ref{HJB.new})}
satisfies all the properties assumed in the derivation of our candidate.

\subsubsection{Solving the HJB equations}

\vspace{0.5em}We now provide conditions under which the heuristic derivation
of the previous section indeed corresponds to a solution of the original
system of HJB equations \textrm{(\ref{HJB bis eq})}. Since we already solved
the problem for $j=1$, we assume here that $j\geq 2$. Let us define 
\begin{equation*}
\overline{v}_{j}:=v_{j}(b_{j}),
\end{equation*}%
and for $x>0$ and $0<\beta \leq $ 1the functions 
\begin{equation*}
\phi _{\beta }(x):=\left( \frac{1+x}{1+(1+\beta )x}\right) ^{\frac{1}{x}-1},%
\text{ }\psi _{\beta }(x):=\frac{\phi _{\beta }(x)-x}{(1-x)\phi _{\beta }(x)}%
.
\end{equation*}

\vspace{0.5em}

\begin{remark}
\label{psi} It is easy to show that the functions $\psi _{\beta }$ can be
extended to continuous functions on $\mathbb{R}_{+}$ which decrease from $1$
to $\frac{1}{2}$ and that for all $x\geq 0$ 
\begin{equation*}
\psi _{1}(x)=\underset{0<\beta \leq 1}{\inf }\psi _{\beta }(x).
\end{equation*}
\end{remark}

\vspace{0.5em} We have the following results.

\begin{proposition}
\label{HJB solution prop} Assume that 
\begin{equation}
\frac{r}{\lambda _{j}}-1\leq \frac{\overline{v}_{j-1}}{b_{j-1}}.
\label{continuation decision}
\end{equation}

\begin{itemize}
\item[$\mathrm{{(i)}}$] The ordinary differential equations (\ref{HJB.new}),
along with (\ref{extrapolation eq}), have unique maximal solutions $v_{j}$
for $j\geq 2$. The functions $v_{j}$ are globally concave, differentiable
everywhere except at $b_{j}$ and twice differentiable everywhere except at $%
b_{j}$ and $b_{j}+b_{j-1}$. The endogenous thresholds $\gamma _{j}\geq
b_{j}+b_{j-1}$ are uniquely determined by 
\begin{equation}
\frac{r}{\lambda _{j}}-1\in \partial v_{j-1}(\gamma _{j}-b_{j}),
\label{subdifferential eq}
\end{equation}%
where $\partial v_{j}(u)$ is the subdifferential of $v_{j}$ at $u$, and
verify 
\begin{equation}
\gamma _{j}\leq b_{j}+\gamma _{j-1}.  \label{oufouf}
\end{equation}

\item[$\mathrm{{(ii)}}$] Assume that the $\lambda _{j}$ verify 
\begin{equation}
\left( v_{j-1}^{\prime }(b_{j-1}^{+})\right) ^{+}\frac{b_{j-1}}{\overline{v}%
_{j-1}}\leq \psi _{1}\left( \frac{r}{\lambda _{j}}\right) ,
\label{hyp.lambda}
\end{equation}%
In that case, the functions $v_{j}$ also verify 
\begin{equation}
v_{j}^{\prime }(u)-v_{j-1}^{\prime }(u-b_{j})\leq 0,\text{ for all }u\geq
b_{j}.  \label{propz}
\end{equation}
\end{itemize}
\end{proposition}

\vspace{0.5em} The proof is rather tedious and relegated to the Appendix.
(See \cite{pages} for interpretations.) The technical condition (\ref%
{hyp.lambda}) restricts the range of admissible values for $\left\{ \lambda
_{j}\right\} _{1\leq j\leq I}$. The left-hand side reflects the kink of the
value function at $b_{j}$, and is less than one by concavity. The condition
is met for sufficiently large values of $\lambda _{j}$, since $\psi _{1}$
converges to 1 near the origin.

\vspace{0.5em} Now since the functions $v_{j}$ constructed in Proposition %
\ref{HJB solution prop} are globally concave, have a derivative which is
greater than $-1$ for $u<\gamma _{j}$ and equal to $-1$ for $u\geq \gamma
_{j}$ and satisfy \textrm{(\ref{propz})}, we can apply the heuristic
arguments of Section \ref{formal derivation} to obtain the following
corollary.

\vspace{0.5em}

\begin{corollary}
\label{coco} Under the assumptions of Proposition \ref{HJB solution prop},
the functions $v_{j}$ constructed in the same Proposition solve the HJB
equations \textrm{(\ref{HJBB})}.
\end{corollary}

\vspace{0.5em} {\noindent \textbf{Proof. }} The only remaining property to
prove is that for $u\geq \gamma _{j}$, we have 
\begin{equation*}
-\left( ru+\lambda _{j}b_{j}\right) v_{j}^{\prime }(u)-j\mu +\lambda
_{j}\left( v_{j}(u)-v_{j-1}(u-b_{j})\right) \geq 0.
\end{equation*}%
We compute 
\begin{align*}
& -\left( ru+\lambda _{j}b_{j}\right) v_{j}^{\prime }(u)-j\mu +\lambda
_{j}\left( v_{j}(u)-v_{j-1}(u-b_{j})\right) \\
& =ru+\lambda _{j}b_{j}-j\mu +\lambda _{j}\left( v_{j}(\gamma _{j})-u+\gamma
_{j}-v_{j-1}(u-b_{j})\right) \\
& =r(u-\gamma _{j})+\lambda _{j}\left( v_{j-1}(\gamma _{j}-b_{j})+\gamma
_{j}-b_{j}-v_{j-1}(u-b_{j})-u+b_{j}\right) \\
& \geq r(u-\gamma _{j})-\lambda _{j}(u-\gamma _{j})\left( 1+\underset{\gamma
_{j}-b_{j}\leq x\leq u-b_{j}}{\sup }v_{j-1}^{\prime }(x)\right) \\
& \geq r(u-\gamma _{j})-\lambda _{j}(u-\gamma _{j})\frac{r}{\lambda _{j}} \\
& =0,
\end{align*}%
where we used the fact that $v_{j-1}$ is concave, that $u\rightarrow
v_{j-1}(u)+u$ is increasing and that $v_{j-1}^{\prime }(\gamma
_{j}-b_{j})\leq \frac{r}{\lambda _{j}}-1$. In particular, this shows that 
\begin{align}
& -\underset{(\theta ,h^{1},h^{2})\in \mathcal{B}^{j}}{\sup }\left\{ \left(
ru+\lambda _{j}\left( h^{1}+(1-\theta )h^{2}\right) \right) v_{j}^{\prime
}(u)+j\mu \right.  \notag  \label{label} \\
& \hspace{6.5em}\left. -\lambda _{j}\left( v_{j}(u)-\theta
v_{j-1}(u-h^{1})\right) \right\} \geq 0,\ u\geq \gamma _{j}.
\end{align}%
\hbox{ }\hfill $\Box $

\subsection{The verification theorem\label{verif}}

\vspace{0.5em}In this subsection, we prove our main result.

\begin{theorem}
\label{main} Let $u_0\leq \gamma_I$ be the reservation utility for the bank.
Then, the optimal contract in $\widetilde{\mathcal{A}}^0(x)$ for the problem 
\textrm{(\ref{problem})} is the contract \ref{contract}.
\end{theorem}

\vspace{1em} We decompose the proof in two parts. First, we show that the
bank can obtain a level of utility $u_0$ and the investors $v_I(u_0)$, for
any $u_0\geq b_I$, with this contract. The second part, reported in
Proposition \ref{2}, shows that for any contract $(D,\theta)$ which makes
the shirking decision $k=0$ incentive-compatible, the utility the investors
can obtain is bounded from above by $v_I(u_0)$, where $u_0$ is the utility
obtained by the bank.

\begin{proposition}
\label{prop.contract} Let the assumptions of Proposition \ref{HJB solution
prop} hold true. For any starting condition $u_{0}>b_{I}$, we define the
process $u_{t}$ as the solution of the following SDE for $j=0$, $\dots $, $%
I-1$ 
\begin{align}
du_{t}& =(ru_{t}-\delta
^{I-N_{t}}(u_{t}))\,dt-h^{1,I-N_{t}}(u_{t})(dN_{t}-\lambda _{I-N_{t}}\,dt) 
\notag  \label{sde} \\[0.5em]
& \hspace{0.9em}-h^{2,I-N_{t}}(u_{t})(dH_{t}-\lambda _{I-N_{t}}(1-\theta
^{I-N_{t}}(u_{t}))\,dt),\text{ }t<\tau .
\end{align}

\vspace{1em} Then, the contract defined by $\left( \delta
^{I-N_{t}}(u_{t}),\theta ^{I-N_{t}}(u_{t})\right) $ is incentive-compatible,
has value $u_{0}$ for the bank and value $v_{I}(u_{0})$ for the investors.
\end{proposition}

\vspace{1em} {\noindent \textbf{Proof. }}First, the drift and volatility in
the SDE \textrm{(\ref{sde})} are clearly Lipschitz. This guarantees the
existence and uniqueness of the solution for all $t$. Moreover, it is also
clear from the definitions of $\delta ^{I-N_{t}}$, $\theta ^{I-N_{t}}$, $%
h^{1,I-N_{t}}$ and $h^{2,I-N_{t}}$ that 
\begin{equation*}
ru_{t}-\delta ^{I-N_{t}}+\lambda _{I-N_{t}}\left( h^{1,I-j}(u_{t})+(1-\theta
^{I-N_{t}}(u_{t}))h^{2,I-N_{t}}(u_{t})\right) \geq 0.
\end{equation*}%
Hence $u_{t}$ remains below $\gamma _{I-N_{t}}$. Moreover, when $N$ jumps,
we have at the time of the jump
\begin{align*}
u_{t} =u_{t^{-}}-h_{t}^{1,I-N_{t^{-}}} =&b_{I-N_{t}}1_{b_{I-N_{t^{-}}}\leq
u_{t^{-}}<b_{I-N_{t}}+b_{I-N_{t^{-}}}}\\
&+(u_{t^{-}}-b_{N_{t^{-}}})1_{b_{I-N_{t}}+b_{I-N_{t^{-}}}\leq u_{t^{-}}\leq \gamma _{I-N_{t^{-}}}}
\\
 \geq& b_{I-N_{t}}.
\end{align*}%
Therefore, we always have $u_{t}\geq b_{I-N_{t}}$ for $t<\tau $. Hence, the
process $u$ is bounded.

\vspace{1em} Moreover, it is clear by construction that this contract makes
the shirking decision $k=0$ incentive-compatible. Indeed, we have after some
calculations for all $j$ 
\begin{equation*}
h^{1,I-N_{t}}(u_{t})+(1-\theta
^{I-N_{t}}(u_{t}))h^{2,I-N_{t}}(u_{t})=b_{I-N_{t}},\text{ }t<\tau ,
\end{equation*}%
which is exactly \textrm{(\ref{IC cond})}.

\vspace{1em} Then, using the equation \textrm{(\ref{PK eq})} for the
continuation utility of the bank obtained with the contract $%
(\delta^{I-N_t}(u_t),\theta^{I-N_t}(u_t))$, we obtain
\begin{align*}
d\left( e^{-rt}(u_{t}^{0}-u_{t})\right) & =e^{-rt}\left(
(h_{t}^{1}-h^{1,I-N_{t}})(dN_{t}-\lambda _{I-N_{t}}dt)\right) \\
& \hspace{0.9em}+e^{-rt}\left(
(h_{t}^{2}-h^{2,I-N_{t}})(dH_{t}-\lambda _{I-N_{t}}(1-\theta
^{I-N_{t}})dt)\right),
\end{align*}
where we suppressed the dependance of $h^{1,I-N_t}$, $h^{2,I-N_t}$ and $\theta^{I-N_t}$ in $u$ for simplicity.

\vspace{0.5em}
Since $h^{1,N_{t}}(u_{t})$ and $h^{2,I-N_{t}}(u_{t})$ are bounded because $%
u_{t}$ is bounded and since $h_{t}^{1}$ and $h_{t}^{2}$ are in the space $%
L^{1}(\mathbb{P})$ by construction, we can take the conditional expectation
above to obtain 
\begin{equation*}
\mathbb{E}_{t}\left[ u_{t+s}^{0}-u_{t+s}\right] =e^{rs}(u_{t}^{0}-u_{t}).
\end{equation*}%
$u^{0}$ remains bounded, because the $\delta ^{j}$ are bounded for all $j$
(recall \textrm{(\ref{continuation utility eq})}) and $u$ is bounded, thus
the left-hand side above must remain bounded. Since $r>0$, letting $s$ go to 
$+\infty $ implies that $u_{t}=u_{t}^{0}$, $\mathbb{P}-a.s.$ and in
particular that the bank overall utility is 
\begin{equation*}
u_{0}^{0}=u_{0}.
\end{equation*}

\vspace{1em} Let us now turn our attention to the investors. Define 
\begin{equation}
G_{t}:=\int_{0}^{t}((I-N_{s})\mu -\delta (u_{s}))ds+v_{I-N_{t}}(u_{t}),
\label{invest}
\end{equation}%
where the $v_{j}$ are those defined in Proposition \ref{HJB solution prop}.
Consider the interval $[\tau _{j}\wedge \tau ,\tau _{j+1}\wedge \tau )$. We
have shown before that $u_{t}$ remains above $b_{I-j}$. But we know by
construction that $v_{I-j}$ is continuous on $[b_{I-j},+\infty )$ and has a
derivative which can be continuously extended on $[b_{I-j},+\infty )$. Hence
we can apply the change of variable formula for locally bounded processes
(see \cite{dellm}, Chapter VI, Section $92$) to obtain for all $t\geq 0$ 
\begin{align}
G_{t}& =v_{I}(u_{0})+\sum_{j=0}^{I-1}\int_{\tau _{j}\wedge t}^{\tau
_{j+1}\wedge t}(I-j)\mu -\delta ^{I-j}(u_{s})+v_{I-j}^{\prime }(u_{s})(
ru_{s}-\delta ^{I-j}(u_{s})) ds  \notag  \label{G} \\
& \hspace{0.9em}+\sum_{j=0}^{I-1}\int_{\tau _{j}\wedge t}^{\tau _{j+1}\wedge
t}\lambda _{I-j}v_{I-j}^{\prime }(u_{s})\left( h^{1,I-j}(u_{s})+(1-\theta
^{I-j}(u_{s}))h^{2,I-j}(u_{s})\right) ds  \notag \\
& \hspace{0.9em}+\sum_{j=0}^{I-1}\sum_{\tau _{j}\wedge t\leq s\leq \tau
_{j+1}\wedge t}v_{I-j}(u_{s})-v_{I-j}(u_{s^{-}}).
\end{align}

\vspace{0.5em} Let us decompose the jumps of $v_{j}$. We have 
\begin{align*}
v_{j}(u_{s})-v_{j}(u_{s^{-}}) =&\Delta N_{s}\left( \left( 1-\Delta
H_{s}\right) v_{j-1}\left( u_{s^{-}}-h^{1,j}(u_{s^{-}})\right) -v_{j}\left(
u_{s^{-}}\right) \right) \\[0.5em]
 =&\Delta N_{s}\left( v_{j-1}\left( u_{s^{-}}-h^{1,j}(u_{s^{-}})\right)
-v_{j}\left( u_{s^{-}}\right) \right)\\
 &-\Delta H_{s}v_{j-1}\left(
u_{s^{-}}-h^{1,j}(u_{s^{-}})\right) ,
\end{align*}%
which implies that 
\begin{align*}
\sum_{s=\tau _{j}\wedge t}^{ \tau _{j+1}\wedge t}%
\scriptstyle v_{I-j}(u_{s})-v_{I-j}(u_{s^{-}})& =\int_{\scriptstyle\tau
_{j}\wedge t}^{\scriptstyle\tau _{j+1}\wedge t}\scriptstyle\left( v_{I-j-1}\left(
u_{s^{-}}-h^{1,I-j}(u_{s^{-}})\right) -v_{I-j}\left( u_{s^{-}}\right)
\right) dN_{s} \\
& \hspace{0.9em}-\int_{\scriptstyle\tau _{j}\wedge t}^{\scriptstyle\tau
_{j+1}\wedge t}v_{I-j-1}\left( u_{s^{-}}-h^{1,I-j}(u_{s^{-}})\right) dH_{s}.
\end{align*}%
From this, we obtain 
\begin{align*}
&G_{t} =v_{I}(u_{0})+\sum_{j=0}^{I-1}\int_{\tau _{j}\wedge t}^{\tau
_{j+1}\wedge t}(I-j)\mu -\delta ^{I-j}(u_{s})+v_{I-j}^{\prime }(u_{s})(
ru_{s}-\delta ^{I-j}(u_{s})) ds \\
& +\sum_{j=0}^{I-1}\int_{\tau _{j}\wedge t}^{\tau _{j+1}\wedge
t}\lambda _{I-j}v_{I-j}^{\prime }(u_{s})\left( h^{1,I-j}(u_{s})+(1-\theta
^{I-j}(u_{s}))h^{2,I-j}(u_{s})\right) ds \\
& +\sum_{j=0}^{I-1}\int_{\tau _{j}\wedge t}^{\tau _{j+1}\wedge
t}\lambda _{I-j}\left( v_{I-j-1}\left( u_{s}-h^{1,I-j}(u_{s})\right)
-v_{I-j}\left( u_{s}\right) \right) ds \\
& -\sum_{j=0}^{I-1}\int_{\tau _{j}\wedge t}^{\tau _{j+1}\wedge
t}\lambda _{I-j}(1-\theta ^{I-j})v_{I-j-1}\left(
u_{s}-h^{1,I-j}(u_{s})\right) ds \\
& +\sum_{j=0}^{I-1}\int_{\tau _{j}\wedge t}^{\tau _{j+1}\wedge
t}\left( v_{I-j-1}\left( u_{s^{-}}-h^{1,I-j}(u_{s^{-}})\right)
-v_{I-j}\left( u_{s^{-}}\right) \right) \left( dN_{s}-\lambda _{I-j}ds\right)
\\
& -\sum_{j=0}^{I-1}\int_{\tau _{j}\wedge t}^{\tau _{j+1}\wedge
t}v_{I-j-1}\left( u_{s^{-}}-h^{1,I-j}(u_{s^{-}})\right) \left(
dH_{s}-\lambda _{I-j}(1-\theta ^{I-j}(u_{s^{-}}))ds\right) .
\end{align*}

\vspace{0.5em} Using the fact that the $v_{j}$ solve the HJB equation \ref%
{HJB.new}, we deduce that%
\begin{align}\label{GG}
\nonumber &G_{t} =v_{I}(u_{0})+\sum_{j=0}^{I-1}\int_{\tau _{j}\wedge t}^{\tau
_{j+1}\wedge t}v_{I-j-1}\left( u_{s^{-}}-h^{1,I-j}(u_{s^{-}})\right) \left( dN_{s}-\lambda _{I-j}ds\right)\\
\nonumber&-\sum_{j=0}^{I-1}\int_{\tau _{j}\wedge t}^{\tau
_{j+1}\wedge t}v_{I-j}\left( u_{s^{-}}\right) \left( dN_{s}-\lambda _{I-j}ds\right)\\
& -\sum_{j=0}^{I-1}\int_{\tau _{j}\wedge t}^{\tau _{j+1}\wedge
t}v_{I-j-1}\left( u_{s^{-}}-h^{1,I-j}(u_{s^{-}})\right) \left(
dH_{s}-\lambda _{I-j}(1-\theta ^{I-j}(u_{s^{-}}))ds\right) .
\end{align}

\vspace{0.5em} Hence, $G$ is a bounded martingale until time $\tau$ (since $%
\delta$ is bounded by definition and $u_t$ and thus the $v_j(u_t)$ are also
bounded) and we have, since $u_\tau=0$ 
\begin{equation*}
\mathbb{E}\left[\int_0^\tau\left((I-N_t)\mu-\delta_t\right)dt\right]=\mathbb{%
E}[G_\tau]=G_0=v_I(u_0),
\end{equation*}
which is the desired result. \hbox{ }\hfill$\Box$

\vspace{0.5em} We now show that $v_I(u_0)$ is an upper bound for the utility
the investor can obtain from any contract which makes the shirking decision $%
k=0$ incentive-compatible.

\begin{proposition}
\label{2} For any contract $(D,\theta)\in\widetilde{\mathcal{A}}^0(u_0)$,
the utility the investors can obtain is bounded from above by $v_I(u_0)$,
where $u_0$ is the utility obtained by the bank.
\end{proposition}

\vspace{0.5em} {\noindent \textbf{Proof. }} We define as in the previous
proof the quantity $G_t$ for an arbitrary contract $(\delta,\theta)$. By
applying the change of variable formula and arguing exactly as before we can
obtain that the drift of $G$ is actually negative, using again \textrm{(\ref%
{HJBB})}. Indeed, we know that for any $(D,\theta,h^1,h^2)\in\widetilde{%
\mathcal{A}}^0(u_0)$, we have from Corollary \ref{coco} and its proof that
for all $j$ 
\begin{equation*}
\left( ru_t+\lambda_j\left(h^1_t+(1-\theta_t)h^2_t\right)\right)v^{\prime
}_j(u_t)+j\mu
-\lambda_j\left(v_j(u_t)-\theta_tv_{j-1}(u_t-h^1_t)\right)\leq0,
\end{equation*}
and we know that 
\begin{equation*}
-(v^{\prime }_j(u_t)+1)dD_t\leq 0,
\end{equation*}
since $D$ is non-decreasing.

\vspace{0.5em} Hence, using again \textrm{(\ref{GG})}, we have 
\begin{align}
G_{t\wedge \tau }& \leq v_{I}(u)+\int_{0}^{\tau \wedge t}\left(
v_{I-N_{s}-1}\left( u_{s^{-}}-h_{s}^{1,I-N_{s}}\right) -v_{I-N_{s}}\left(
u_{s^{-}}\right) \right)\scriptstyle \left( dN_{s}-\lambda _{I-N_{s}}ds\right)  \notag
\label{GGG} \\
& \hspace{0.9em}-\int_{0}^{\tau \wedge t}v_{I-N_{s}-1}\left(
u_{s^{-}}-h_{s}^{1,I-N_{s}}\right) \left( dH_{s}-\lambda _{I-N_{s}}(1-\theta
_{s}^{I-N_{s}})ds\right) .
\end{align}

\vspace{0.5em}Now we have 
\begin{align*}
& \mathbb{E}\left[ \int_{0}^{\tau \wedge t}\left\vert v_{I-N_{s}-1}\left(
u_{s}-h_{s}^{1,I-N_{s}}\right) -v_{I-N_{s}}\left( u_{s}\right) \right\vert ds%
\right] \\
& \leq \mathbb{E}\left[ \int_{0}^{\tau \wedge t}\left\vert
v_{I-N_{s}-1}\left( u_{s}-h_{s}^{1,I-N_{s}}\right) -v_{I-N_{s}-1}\left(
u_{s}-b_{I-N_{s}}\right) \right\vert ds\right] \\
& \hspace{0.9em}+\mathbb{E}\left[ \int_{0}^{\tau \wedge t}\left\vert
v_{I-N_{s}-1}\left( u_{s}-b_{I-N_{s}}\right) -v_{I-N_{s}}\left( u_{s}\right)
\right\vert ds\right] \\
&
\end{align*}

\vspace{0.5em} Then, from \textrm{(\ref{propz})}, we know that for all $j$
the function $u\longrightarrow v_{j}(u)-v_{j-1}(u-b_{j})$ is decreasing.
Moreover, for $u$ large enough (namely $u\geq \gamma _{j}\vee (\gamma
_{j-1}+b_{j})$) we have 
\begin{equation*}
v_{j}(u)-v_{j-1}(u-b_{j})=v_{j}(\gamma _{j})+\gamma _{j}-v_{j-1}(\gamma
_{j-1})+\gamma _{j-1}-b_{j},
\end{equation*}%
which implies that for all $j$ the function $u\longrightarrow
v_{j}(u)-v_{j-1}(u-b_{j})$ is bounded. Moreover, we have 
\begin{align*}
& \mathbb{E}\left[ \int_{0}^{\tau \wedge t}\left\vert v_{I-N_{s}-1}\left(
u_{s}-h^{1,I-N_{s}}\right) -v_{I-N_{s}-1}\left( u_{s}-b_{I-N_{s}}\right)
\right\vert ds\right] \\
& \leq \mathbb{E}\left[ \int_{0}^{\tau \wedge t}\left\vert
h_{s}^{1,I-N_{s}}-b_{I-N_{s}}\right\vert \underset{b_{I-N_{s}}<u\leq \gamma
_{I-N_{s}}}{\sup }\left\vert v_{I-N_{s}}^{\prime }(u)\right\vert ds\right] \\
& \leq C\left( 1+\mathbb{E}\left[ \int_{0}^{\tau \wedge t}\left\vert
u_{s}\right\vert ds\right] \right) \\
& \leq C\left( 1+\mathbb{E}\left[ \int_{0}^{\tau \wedge t}ue^{(r+2\lambda
)s}ds\right] \right) <+\infty ,
\end{align*}%
where $\lambda :=\underset{1\leq j\leq I}{\sup }\lambda _{j}$, and where we
used successively the fact that the derivative of the $v_{j}$ can be
extended to a continuous function on $[b_{j},\gamma _{j}]$ which is
therefore bounded on that compact, then the fact that by the limited
liability condition \textrm{(\ref{limited liability eq})} we have $%
h_{t}^{1}\leq u_{t}$, and finally that conditionally on the fact that there
are $j$ loans left in the pool, the drift of $u_{t}$ as given by \textrm{(%
\ref{PK eq})} is%
\begin{align*}
ru_{t}+\lambda _{j}\left( h_{t}^{1}+(1-\theta _{t})h_{t}^{2}\right) -\delta
_{t}& \leq ru_{t}+\lambda _{j}\left( h_{t}^{1}+(1-\theta
_{t})(u_{t}-h_{t}^{1})\right) \\
& \leq ru_{t}+\lambda _{j}\left( u_{t}-b_{j-1}+(1-\theta _{t})u_{t})\right)
\\
& \leq u_{t}(r+2\lambda _{j}),
\end{align*}%
where we used the fact that $h$, $b_{j}$ and $\lambda _{j}$ are positive.
Hence, $u_{t}$ increases at a rate lower than $r+2\lambda $.

\vspace{0.5em} Similarly, we have 
\begin{align*}
&\mathbb{E}\left[\int_{0}^{\tau\wedge
t}\left|v_{I-N_s-1}\left(u_{s^-}-h_s^{1,I-N_s}\right)\right|ds\right] \\
&\leq \mathbb{E}\left[\int_{0}^{\tau\wedge t}\left|h^{2,I-N_s}_s\right|%
\underset{b_{I-N_s}<u\leq\gamma_{I-N_s}}{\sup}\left| v^{\prime
}_{I-N_s-1}(u)\right|ds\right] \\
&\leq \mathbb{E}\left[\int_{0}^{\tau\wedge t}\left|u_s\right|\underset{%
b_{I-N_s}<u\leq\gamma_{I-N_s}}{\sup}\left| v^{\prime }_{I-N_s-1}(u)\right|ds%
\right]<+\infty.
\end{align*}

\vspace{0.5em} Taking expectations in \textrm{(\ref{GGG})}, we therefore
obtain 
\begin{align}  \label{GGGG}
\nonumber v_I(u_0)\geq& \mathbb{E}\left[\int_0^\tau\left((I-N_s)\mu-\delta_s\right)ds%
\right]\\
\nonumber&+\mathbb{E}\left[1_{t<\tau}\left(\int_t^\tau\left(\delta_s-(I-N_s)\mu%
\right)ds+v_{I-N_t}(u_t)\right)\right]  \\
\nonumber=&\mathbb{E}\left[\int_0^\tau\left((I-N_s)\mu-\delta_s\right)ds\right]\\
&\nonumber+%
\mathbb{E}\left[1_{t<\tau}\mathbb{E}_t\left[\int_t^\tau\left(%
\delta_s-(I-N_s)\mu\right)ds+v_{I-N_t}(u_t)\right]\right] \\
\nonumber=&\mathbb{E}\left[\int_0^\tau\left((I-N_s)\mu-\delta_s\right)ds\right]\\
&\nonumber+%
\mathbb{E}\left[1_{t<\tau}\left(u_t+v_{I-N_t}(u_t)-\mathbb{E}_t\left[%
\int_t^\tau\left(I-N_s\right)\mu ds\right]\right)\right]  \\
\geq& \mathbb{E}\left[\int_0^\tau\left((I-N_s)\mu-\delta_s\right)ds\right]+%
\mathbb{E}\left[1_{t<\tau}\left(-I\mu\tau+u_t+v_{I-N_t}(u_t)\right)\right].
\end{align}

\vspace{0.5em} Then, we know that for all $j$ the function $u\longrightarrow
u+v_j(u)$ is increasing before $\gamma_j$ and is constant for $u\geq\gamma_j$%
. It is therefore bounded and we have 
\begin{equation*}
\left|-I\mu\tau+u_t+v_{I-N_t}(u_t)\right|\leq I\mu\tau +\underset{1\leq
j\leq I}{\sup}\left|\gamma_j+v_j(\gamma_j)\right|\leq C(1+\tau),
\end{equation*}
for some positive constant $C$. This quantity being integrable, we can apply
the dominated convergence theorem in \textrm{(\ref{GGGG})} and let $t$ go to 
$+\infty$ to obtain 
\begin{equation*}
v_I(u_0)\geq \mathbb{E}\left[\int_0^\tau\left((I-N_s)\mu-\delta_s\right)ds%
\right],
\end{equation*}
which is the desired result. \hbox{ }\hfill$\Box$

\subsection{Numerical results}

\vspace{0.5em}\label{sec.num} In this section we present some numerical
results to illustrate our main properties. Following the empirical estimates
of \cite{laurent}, we choose to work with a pool of $I=30$ loans with:

\begin{equation*}
\begin{tabular}{|c|c|}
\hline
$\mu $ & $0.06$ \\ \hline
$r$ & $0.02$ \\ \hline
$B$ & $0.002$ \\ \hline
$\varepsilon $ & $0.25$ \\ \hline
$\left( \alpha _{j}\right) _{1\leq j\leq 20}$ & $0.055$ \\ \hline
$\left( \alpha _{j}\right) _{21\leq j\leq 26}$ & $0.05$ \\ \hline
$\left( \alpha _{j}\right) _{27\leq j\leq 30}$ & $0.044.$ \\ \hline
\end{tabular}%
\end{equation*}

\vspace{0.5em}Assumptions \ref{assump.mu} to \ref{decreasing assumption} are
satisfied. The values assigned to $\alpha _{j}$ and $\epsilon $ are
consistent with the literature. The former are assumed piecewise constant to
model a surge in the default intensity after a certain fraction of the pool
has defaulted. The bank's discount rate and the yield of the loans are taken
close to what could be deemed standard in financial markets.

\begin{figure}[htb]\label{fig.v}
\begin{center}
\centering
\includegraphics[width=30em,height=20em]{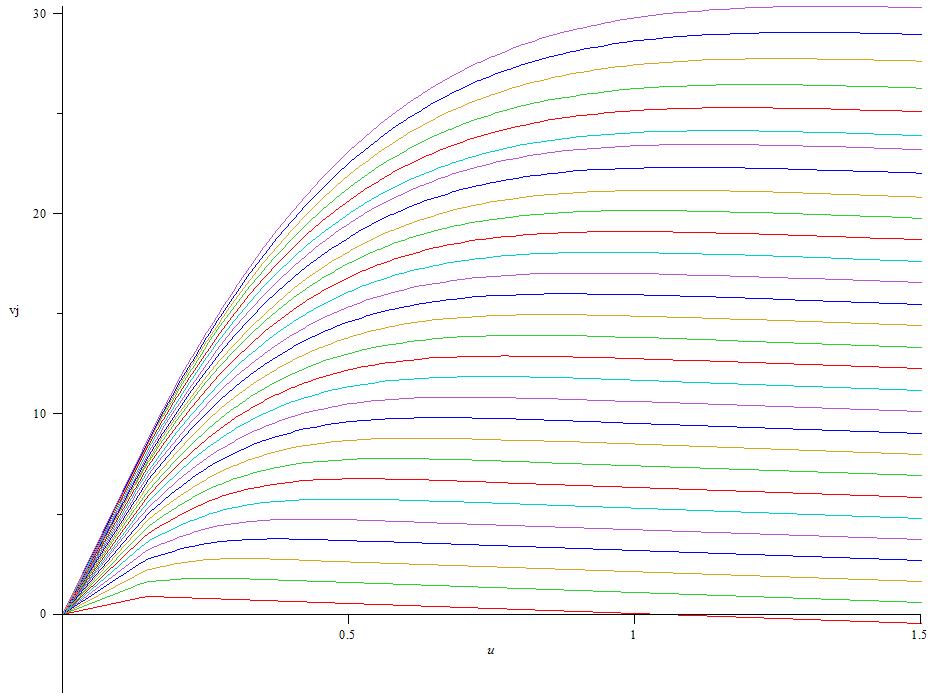}
\caption{Functions $v_j(u)$ for $j=2..30$.}
\end{center}
\end{figure}

\begin{figure}[htb]\label{fig.g}
\begin{center}
\centering
\includegraphics[width=30em,height=20em]{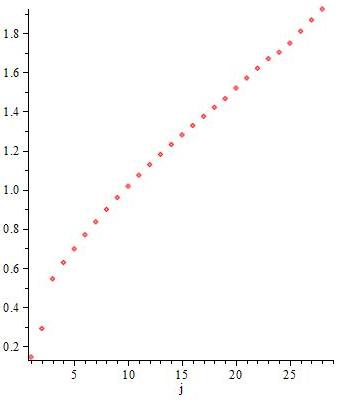}
\caption{Values of $\gamma_j$ for $j=1..30$.}
\end{center}
\end{figure}

\vspace{0.5em}Using the fact that the $v_{j}$ have a semi-explicit form, we
use numerical integrations techniques to obtain the functions $v_{j}$ for $%
j=1$, $\dots $, $30$. Condition \textrm{(\ref{hyp.lambda})} of Proposition %
\ref{HJB solution prop} is always verified. As shown in Figures \ref{fig.v}
and \ref{fig.g}, both the value functions $v_{j}$ and the thresholds $%
\gamma _{j}$ appear to be increasing with $j$. Note that the contract
between the bank and the investors generates a positive social surplus,
given by 
\begin{equation*}
v_{30}(\gamma _{30})+\gamma _{30}-30=1.86.
\end{equation*}%
With competitive investors, the full surplus is extracted by the bank in the
form of expected profits. We have $v_{30}(\gamma _{30})=29.88$ in this
numerical example. Since investors must break even, the bank does not have
the wherewithal to go for the project on its own. The capital that it has to
invest corresponds to roughly $0.4\%$ of the total amount. By construction,
it is the stake that it is willing to invest at time 0 in order to maximize
its profits from the pool. Note that from (\ref{first-best}) the surplus
available under the first-best is $S=4.48$. About 58\% of value is lost due
to the agency problem.

%

\vspace{0.5em}

%
%

\section{What happens when $r=0$?\label{rr}}

\vspace{0.5em}In this section, we relax the assumption that the bank is
impatient. A positive discount rate implies that there are gains from
trades, as the bank is eager to sell claims on future cash flows to more
patient investors. This is precisely what creates the moral hazard problem,
as gains from trades can be undermined by high default rates when the bank
shirks. We will see that, in contrast, the first-best is attained when $r=0$%
. Proofs are quite similar and only sketched. The analog of Proposition \ref%
{HJB solution prop} is as follows.

\begin{proposition}
\label{HJB solution prop2} Assume that $r=0$.\newline
\vspace{-0.2in}

\begin{itemize}
\item[$\mathrm{{(i)}}$] The ordinary differential equations (\ref{HJB.new}),
along with (\ref{extrapolation eq}), have unique maximal solutions $v_{j}$
for $j\geq 1$. The functions $v_{j}$ are globally concave, differentiable
everywhere except at $b_{j}$ and twice differentiable everywhere except at $%
b_{j}$ and $b_{j}+b_{j-1}$. The endogenous thresholds $\gamma _{j}$ are
uniquely determined by 
\begin{equation}
\gamma _{j}=\frac{jB}{\epsilon \overline{\alpha }_{j}}.
\label{subdifferential eq2}
\end{equation}

\item[$\mathrm{{(ii)}}$] We also have 
\begin{equation}
v_{j}^{\prime }(u)-v_{j-1}^{\prime }(u-b_{j})\leq 0,\text{ for all }u\geq
b_{j}.  \label{propz2}
\end{equation}
\end{itemize}
\end{proposition}

\vspace{0.5em}{\noindent \textbf{Proof. }} \textbf{$\mathrm{(i)}$} When $r=0$%
, the solution of \textrm{(\ref{HJB.new})} for a given $\gamma \geq b_{j}$
is 
\begin{align}
v_{j}(u)& =\frac{j\mu }{\lambda _{j}}+e^{\frac{u-\gamma }{b_{j}}}\left(
v_{j-1}(\gamma -b_{j})-b_{j}\right) +\int_{u}^{\gamma }\frac{e^{\frac{u-x}{%
b_{j}}}}{b_{j}}v_{j-1}(x-b_{j})dx,\text{ }b_{j}<u\leq \gamma
\label{specific vj eq} \\
v_{j}(u)& =\gamma -u+v_{j}(\gamma ),\text{ }u>\gamma .  \notag
\end{align}%
Using the same arguments as in the proof of Proposition \ref{HJB solution
prop}, it is easily proved that the choice of $\gamma $ leading to the
maximum solution is%
\begin{equation*}
\gamma _{j}=\gamma _{j-1}+b_{j}.
\end{equation*}%
Reasoning by induction, we can then prove similarly that the functions $%
v_{j} $ verify all the desired properties. Moreover, since $\gamma
_{1}=b_{1} $, we obtain that 
\begin{equation*}
\gamma _{j}=\sum_{i=1}^{j}b_{i}=\frac{jB}{\epsilon \overline{\alpha }_{j}}.
\end{equation*}

\vspace{0.5em} \textbf{$\mathrm{(ii)}$} We can prove that

\begin{align*}
v_{j}^{\prime }(u)& =\int_{u}^{\gamma _{j}}\frac{e^{\frac{u-x}{b_{j}}}}{b_{j}%
}\frac{dv_{j-1}}{du}(x-b_{j})dx-e^{\frac{u-\gamma _{j}}{b_{j}}},\text{ }%
b_{j}<u\leq \gamma _{j} \\
\frac{dv_{j}}{du}(u)& =-1,\text{ }u>\gamma _{j}.
\end{align*}%
By the concavity of $v_{j-1}$, this implies that for $b_{j}<u\leq \gamma
_{j} $ 
\begin{equation*}
v_{j}^{\prime }(u)-v_{j-1}^{\prime }(u-b_{j})\leq -e^{\frac{u-\gamma _{j}}{%
b_{j}}}\left( v_{j-1}^{\prime }(u-b_{j})+1\right) \leq 0.
\end{equation*}%
Since \textrm{(\ref{propz2})} is clear when $u>\gamma _{j}$, this proves $%
\mathrm{(ii)}$. \hbox{ }\hfill $\Box $

\vspace{0.5em} Thanks to Proposition \ref{HJB solution prop2}, we have a
concave solution of the HJB equation, then using the same techniques as in
the case $r>0$, we can verify that the optimal contract is given by

\begin{contract}
\label{contract2}When $r=0$, the optimal contract can be described as
follows:

\begin{itemize}
\item[$\mathrm{{(i)}}$] If, at some point, $u_{t}=\gamma _{j}$, there is no
longer any stochastic liquidation. Fees are paid continuously to the bank,
until extinction of the pool, at the rate $\delta _{t}^{i}=iB/\epsilon $ for
all $i\leq j$.

\item[$\mathrm{{(ii)}}$] Otherwise, the policy is the same as in contract %
\ref{contract} (with $r=0$).
\end{itemize}
\end{contract}

\vspace{0.5em}When the bank starts with reservation utility $\gamma _{I}$
(which is the market outcome when investors are competitive), payments are
never suspended since the bank always operates at the thresholds $\gamma
_{j} $ where payments are made. Hence, there can be no stochastic
liquidation. More specifically, we have 
\begin{equation*}
v_{j}(\gamma _{j})=\frac{j\mu }{\lambda _{j}}-b_{j}+v_{j-1}(\gamma _{j-1}),
\end{equation*}%
implying that the social value of the contract is 
\begin{equation*}
\gamma _{I}+v_{I}(\gamma _{I})=\gamma _{I}+\frac{I}{\overline{\alpha }_{I}}%
\left( \mu -\frac{B}{\epsilon }\right) =\frac{I\mu }{\overline{\alpha }_{I}}.
\end{equation*}%
But according to (\ref{first-best}) this is the social value attained in the
first-best. Hence, when the bank is infinitely patient, the first-best is
attained. For each $j$, the bank captures the maximum value of its rent, $%
\sum_{\left\{ i\leq j\right\} }b_{i}$, but this is not socially costly since
there is no loss arising from any \textquotedblleft
rent-preserving\textquotedblright\ fee.

\vspace{0.5em} Note that, to make the problem interesting, we have assumed
that investments are not self-financing. Otherwise, the bank would be free
to invest arbitrarily large amounts and there would be no demand for
investors' liquidity. This means that 
\begin{equation}
v_{I}(\gamma _{I})<I.  \label{surplus}
\end{equation}%
In the general case, such a condition is difficult to work out, but when $%
r=0 $ it is easily shown to be%
\begin{equation*}
\mu -\frac{B}{\varepsilon }<\overline{\alpha }_{I}.
\end{equation*}%
This yields a lower bound for $B$. The moral hazard problem has to be severe
enough that there is a funding problem. In the general case, we expect that
the equivalent of \textrm{(\ref{surplus})} is going to hold for some lower
bound on $B$ (which will depend on $r$).

\begin{appendix}
\section{Appendix}

\proof[Proof of Proposition \ref{prop.un}]
In this particular case, Problem (\ref{Principle
agent pb}) becomes%
\begin{align}
& v_{1}(u)=\sup_{D}\mathbb{E}^{\mathbb{P}}\left[ \int_{0}^{\tau }\mu
\,dt-dD_{t}\right]  \label{eq.hjbloan} \\
\text{subject to }& \mathbb{E}^{\mathbb{P}}\left[ \int_{0}^{\tau
}e^{-rt}\,dD_{t}\right] \geq u  \notag \\
& u_{t}\geq b_{1},\text{ for all }t<\tau .  \notag
\end{align}

Consider first the subproblem derived from (\ref{eq.hjbloan}) by abstracting
from the initial payment $D_{0}$ and ignoring the incentive compatibility
constraint $u_{t}\geq b_{1}$:%
\begin{align*}
& \widetilde{v}_{1}(u)=\sup_{D}\mathbb{E}^{\mathbb{P}}\left[ \int_{0+}^{\tau
}\mu \,dt-dD_{t}\right] \\
\text{subject to }& \mathbb{E}^{\mathbb{P}}\left[ \int_{0+}^{\tau
}e^{-rt}\,dD_{t}\right] \geq u.
\end{align*}%
The constraint can be written equivalently 
\begin{equation*}
\mathbb{E}^{\mathbb{P}}\left[ \int_{0^{+}}^{\tau }e^{-rt}\left(
dD_{t}-(r+\lambda _{1})u\,dt\right) \right] \geq 0.
\end{equation*}%
The corresponding Lagrangian is 
\begin{equation*}
\mathcal{L}_{t}=\mu \,dt-dD_{t}+\nu _{t}e^{-rt}\left( dD_{t}-u(r+\lambda
_{1})\,dt\right) ,
\end{equation*}%
where $\nu _{t}$ is the Lagrange multiplier at time $t$. Optimizing with
respect to $D$, we get $\nu _{t}=e^{rt}$ and the complementary slackness
conditions imply that the dividend process is absolutely continuous and
constant, namely $dD_{t}=\delta _{t}\,dt$, with $\delta _{t}=(r+\lambda
_{1})u$. Since the process $D$ thus obtained is clearly admissible, this
yields $\widetilde{v}_{1}(u)=\left( \mu -(r+\lambda _{1})u\right) /\lambda
_{1}$.

\vspace{0.5em} Turning now to (\ref{eq.hjbloan}), but still ignoring the
incentive compatibility constraint, we have 
\begin{equation*}
v_{1}(u)=\sup_{D_{0}}-D_{0}+\widetilde{v}_{1}(u-D_{0}),
\end{equation*}
which is increasing in $D_{0}$ when $r>0$. Since $u_{0}=u-D_{0}$ from the
bank's promise-keeping constraint (\ref{PK eq}), the highest initial payment
consistent with the incentive compatibility constraint at time 0 is $%
D_{0}=u-b_{1}$. This yields%
\begin{eqnarray*}
v_{1}(u) &=&b_{1}-u+\widetilde{v}_{1}(b_{1}) \\
&=&b_{1}-u+\overline{v}_{1},
\end{eqnarray*}%
where $\overline{v}_{1}$ is defined as in the Proposition. Finally, one
verifies that $\delta _{t}=b_{1}(r+\lambda _{1})$ yields $u_{t}=b_{1}$ on $%
\left[ 0,\tau \right) $, so that the incentive compatibility condition binds
at all times before default, as desired. 
\qed

\vspace{0.5em}
\proof[Proof of Proposition \ref{HJB solution prop}$\rm{(i)}$]
We will show the result by induction. 
\begin{itemize}
\item Initialization with $j=2$
\end{itemize}

\vspace{0.5em}
The solution of the ODE \reff{HJB.new} for $j=2$ and a given fixed value of $\gamma\geq b_2$ can be easily calculated and is given by

\begin{align*}
\widetilde{v}_2(u,\gamma)&:=(ru+\lambda_2b_2)^{\frac{\lambda_2}{r}}\int_u^{\gamma}\frac{2\mu+\lambda_2v_1(x-b_2)}{(rx+\lambda_2b_2)^{\frac{\lambda_2}{r}+1}}dx\\
&\hspace{0.9em}+\left(v_1(\gamma-b_2)+\frac{2\mu-(r\gamma+\lambda_2b_2)}{\lambda_2}\right)\left(\frac{ru+\lambda_2b_2}{r\gamma+\lambda_2 b_2}\right)^{\frac{\lambda_2}{r}},\text{ }b_2<u\leq\gamma,
\end{align*}
and $\widetilde{v}_2(u,\gamma)=\gamma-u+v_2(\gamma)$ for $u>\gamma$.

\vspace{0.5em}
Now since we have shown that $v_1$ is everywhere twice differentiable except at $b_1$, we have for every $\gamma\neq b_1+b_2$ and every $b_2<u\leq\gamma$
$$\frac{\partial \widetilde{v}_2}{\partial\gamma}(u,\gamma)=\left(v'_1(\gamma-b_2)+1-\frac{r}{\lambda_2}\right)\left(\left(\frac{ru+\lambda_2b_2}{r\gamma+\lambda_2b_2}\right)^{\frac{\lambda_2}{r}}1_{u\leq \gamma}+1_{u>\gamma}\right).$$

\vspace{0.5em}
Thus, the above expression always has the sign of $v'_1(\gamma-b_2)+1-\frac{r}{\lambda_2}$, that is to say that it is positive for $\gamma<b_1+b_2$ 
and negative for $\gamma>b_1+b_2$. Hence, we clearly have for all $b_2< u$
$$\underset{\gamma\geq b_2}{\sup}\widetilde{v}_2(u,\gamma)=\widetilde{v}_2(u,b_1+b_2),$$
which means that the maximal solution of \reff{HJB.new} for $j=2$ corresponds to the choice $\gamma_2=b_1+b_2$, which also happens to correspond to the unique solution of 
$$\frac{r}{\lambda _{2}}-1\in \partial v_{1}(\gamma _{2}-b_1).$$

\vspace{0.5em}
Then, after some calculations, we obtain that for all $b_2<u<b_1+b_2$
$$v''_2(u)=-\left(\lambda_2-r+\lambda_2\frac{\overline{v}_1}{b_1}\right)\frac{ \left(ru+\lambda_2b_2\right)^{\frac{\lambda_2}{r}-1}}{\left(r(b_1+b_2)+\lambda_2b_2\right)^{\frac{\lambda_2}{r}}}\leq0,$$
because of \reff{continuation decision}.

\vspace{0.5em}
Hence, since $v_2$, is linear on $[b_1+b_2,+\infty)$ and is differentiable at $b_1+b_2$, it is concave on $(b_2,+\infty)$. Now if we consider the linear extrapolation of $v_2$ over $[0,b_1]$ by \reff{extrapolation eq}, we just need to verify that the left-derivative of $v_2$ at $b_2$ is less than its right-derivative to obtain the concavity of $v_2$ over $[0,+\infty]$. Taking the limit for $u\downarrow b_2$ in the equation \reff{HJB.new}, we obtain
$$v'_2(b_2^+)=\frac{\lambda_2\overline{v}_2-2\mu}{b_2(r+\lambda_2)}.$$
This implies that
$$v'_2(b_2^-)-v'_2(b_2^+)=\frac{2\mu}{b_2\lambda_2}+v'_2(b_2^+)\frac{r}{\lambda_2}\geq \frac{\mu\epsilon}{B}
-\frac{r}{\lambda_2}.$$

\vspace{0.5em}
Now recall Assumption \ref{assump.alpha}, which implies that
$$\frac{r}{\lambda_j}<\frac{r}{\overline{\alpha}_j}\leq\frac{\mu\epsilon-B}{B}\frac{\epsilon}{1+\epsilon}<\frac{\mu\epsilon}{B}$$
for any $j$ so that $v'_2(b_2^-)-v'_2(b_2^+)\geq 0$.

\begin{itemize}
\item Heredity : $j\geq 3$
\end{itemize}

Let us now suppose that the maximal solution of \reff{HJB.new} $v_{j-1}$ has been constructed for some $j\geq 3$, that it is globally concave on $[0,+\infty)$, everywhere differentiable except at $b_{j-1}$, everywhere twice differentiable except at $b_{j-1}$ and $b_{j-1}+b_{j-2}$, and that the corresponding $\gamma_{j-1}\geq b_{j-1}+b_{j-2}$. Let us now construct the maximal solution corresponding to $j$. Exactly as in the case $j=2$, the solution of the ODE \reff{HJB.new} and a given fixed value of $\gamma\geq b_j$ can be easily calculated and is given by

\begin{align*}
\widetilde{v}_j(u,\gamma):=&(ru+\lambda_jb_j)^{\frac{\lambda_j}{r}}\int_u^{\gamma}\frac{j\mu+\lambda_jv_{j-1}(x-b_j)}{(rx+\lambda_jb_j)^{\frac{\lambda_j}{r}+1}}dx\\
&+\left(v_{j-1}(\gamma-b_j)+\frac{j\mu-(r\gamma+\lambda_jb_j)}{\lambda_j}\right)\left(\frac{ru+\lambda_jb_j}{r\gamma+\lambda_j b_j}\right)^{\frac{\lambda_j}{r}},\text{ }b_j<u\leq\gamma,
\end{align*}
and $\widetilde{v}_j(u,\gamma)=\gamma-u+v_j(\gamma)$ for $u>\gamma$.

Note also that from \reff{HJB.new} it is clear that $v_j$ is differentiable everywhere except at $b_j$,  and twice differentiable everywhere except at $b_j$ and $b_j+b_{j-1}$.

\vspace{0.5em}
Now since we assumed that $v_{j-1}$ is everywhere differentiable except at $b_{j-1}$, we have for every $\gamma\neq b_{j-1}+b_j$ and every $b_j<u\leq\gamma$
$$\frac{\partial \widetilde{v}_j}{\partial\gamma}(u,\gamma)=\left(v'_{j-1}(\gamma-b_j)+1-\frac{r}{\lambda_j}\right)\left(\left(\frac{ru+\lambda_jb_j}{r\gamma+\lambda_jb_j}\right)^{\frac{\lambda_j}{r}}1_{u\leq \gamma}+1_{u>\gamma}\right).$$

\vspace{0.5em}
Thus, since $v_{j-1}$ is concave and its derivative non-increasing, we can conclude as in the case $j=2$ that the maximal solution is uniquely determined by the choice $\gamma_j$ which corresponds to the solution of
$$\frac{r}{\lambda _{j}}-1\in \partial v_{j-1}(\gamma _{j}-b_j).$$

\vspace{0.5em}
More precisely, using \reff{continuation decision}, we have only two cases. Either,
$$v'_{j-1}(b_{j-1}^+)\leq \frac{r}{\lambda_j}-1\leq\frac{\overline{v}_{j-1}}{b_{j-1}},$$ 
and $\gamma_j=b_{j-1}+b_j$, or 
$$\frac{r}{\lambda_j}-1<v'_{j-1}(b_{j-1}^+),$$
and $b_{j-1}+b_j<\gamma_j\leq\gamma_{j-1}+b_j$.

\vspace{0.5em}
Let us now study the concavity. We can differentiate twice the equation \reff{HJB.new} on $(b_j,b_j+b_{j-1})$ since $v_{j-1}(u-b_j)$ is linear and thus twice differentiable on this open interval. We then obtain easily

\begin{align}\label{eq..}
v''_j(u)=v''_j((b_j+b_{j-1})^-)\left(\frac{ru+\lambda_jb_j}{r(b_j+b_{j-1})+\lambda_jb_j}\right)^{\frac{\lambda_j}{r}-2}, \text{ }b_j<u<b_j+b_{j-1}.
\end{align}

\vspace{0.5em}
There are then two cases. If $\gamma_j=b_j+b_{j-1}$, differentiating once \reff{HJB.new} and then taking the limit $u\uparrow b_j+b_{j-1}$, we get
$$(r(b_j+b_{j-1})+\lambda_jb_j)v''_j((b_j+b_{j-1})^-)=\lambda_j\left(\frac{r}{\lambda_j}-1-\frac{\overline{v}_{j-1}}{b_{j-1}}\right)\leq 0.$$

\vspace{0.5em}
Since $v''_j(u)=0$ for $u>b_j+b_{j_1}$, we have proved the concavity on $(b_j,+\infty)$.

\vspace{0.5em}
Now if $\gamma_j>b_j+b_{j-1}$, differentiating once \reff{HJB.new} and taking limits on both sides of $b_j+b_{j-1}$, we obtain

\begin{equation}\label{eq.vpoint2}
v''_j((b_j+b_{j-1})^+)-v''_j((b_j+b_{j-1})^-)=\frac{\lambda_j}{r(b_j+b_{j-1})+\lambda_jb_j}\left(\frac{\overline{v}_{j-1}}{b_{j-1}}-v'_{j-1}(b_{j-1}^+)\right),
\end{equation}
where the right-hand side is positive by the concavity of $v_{j-1}$.

\vspace{0.5em}
Next, we differentiate twice \reff{HJB.new} on $(b_j+b_{j-1},\gamma_j]$. We obtain easily

\begin{equation}\label{v2point}
v''_j(u)=\lambda_j(ru+\lambda_jb_j)^{\frac{\lambda_j}{r}-2}\int_u^{\gamma_j}\frac{v''_{j-1}(x-b_j)}{(ru+\lambda_jb_j)^{\frac{\lambda_j}{r}-1}}dx.
\end{equation}

\vspace{0.5em}
Note that we should normally distinguish between the cases $b_j+b_{j-1}+b_{j-2}\leq\gamma_j$ or not, since $v_{j-1}$ is not twice differentiable at $b_{j-1}+b_{j-2}$. However, since we know that $v_j$ is twice differentiable at $b_j+b_{j-1}+b_{j-2}$, this actually does not change the result. Since $v_{j-1}$ is concave, \reff{v2point} implies that $v_j$ is concave on $(b_j+b_{j-1},+\infty)$. Then with \reff{eq.vpoint2} we obtain that the left second derivative of $v_j$ at $b_j+b_{j-1}$ is negative, which, thanks to \reff{eq..} shows finally the concavity on $(b_j,+\infty)$.

\vspace{1.2em}
Finally, it remains to show that $v'_j(b_j^+)\leq \frac{\overline{v}_j}{b_j}$. We take the limit for $u\downarrow b_j$ in the equation \reff{HJB.new}, we obtain

$$v'_j(b_j^+)=\frac{\lambda_j\overline{v}_j-j\mu}{b_j(r+\lambda_j)}.$$

Since $v'_j\geq -1$, this implies that
$$v'_j(b_j^-)-v'_j(b_j^+)=\frac{j\mu}{b_j\lambda_j}+v'_j(b_j^+)\frac{r}{\lambda_j}\geq \frac{\mu\epsilon}{B}-\frac{r}{\lambda_j},$$
which has already been shown to be positive under Assumption \ref{assump.alpha}. Hence $v_j$ is concave on $[0,+\infty)$. 
\qed

\vspace{1em}
\proof[Proof of Proposition \ref{HJB solution prop}$\rm{(ii)}$]
First of all, by the properties of the function $\psi_1$ recalled in Remark \ref{psi}, it is clear that we can always find a $\lambda_j$ such that \reff{hyp.lambda} is satisfied. Then, if for a fixed $j\geq 2$ we have $v'_{j-1}(b_{j-1}^+)\leq 0$, by differentiating \reff{HJB.new}, we immediately have for $u>b_j$ and $u\neq b_j+b_{j-1}$
\begin{equation}\label{machinchose}
\lambda_j\left(v'_j(u)-v'_{j-1}(u-b_j)\right)=(ru+\lambda_jb_j)v''_j(u)+rv'_j(u).
\end{equation}

Since we have proved in $\rm{(i)}$ that the $v_j$ are concave, it is clear that if $v'_{j-1}(b_{j-1}^+)\leq 0$, the right-hand side above is negative. Then by left and right continuity of $v'_{j-1}$ at $b_{j-1}$, the result extends to $u=b_j+b_{j-1}$. Hence the desired property \reff{propz}. In particular, this proves the result for $j=2$ since $v'_{1}(b_{1}^+)=-1$.

\vspace{0.5em}
Note also that the property \reff{propz} clearly holds for $v_j$ when $u>\gamma_j$. Indeed, we have 
$$v'_{j}=-1$$ and we know that the derivative of $v_{j-1}$ is always greater than $-1$.

\vspace{0.5em}
Let us now show the rest of the result by induction. Since \reff{propz} is true for $j=2$, let us fix a $j\geq 3$ and assume that 
\begin{equation}\label{hdr}
v'_{j-1}(u)-v'_{j-2}(u-b_{j-1})\leq 0,\text{ }u>b_{j-1}.
\end{equation}

\vspace{0.5em}
Now if $v'_{j-1}(b_{j-1}^+)\leq 0$, we already know that the property \ref{propz} is true for $v_j$, so we will assume that $v'_{j-1}(b_{j-1}^+)> 0$. Moreover, by our remark above, we know that \reff{propz} holds true for $v_j$ when $u>\gamma_j$. Let us then first prove that \reff{propz} for $v_j$ when $u>b_j+b_{j-1}$. If $\gamma_j=b_j+b_{j-1}$, there is nothing to do. Otherwise, we have using successively \reff{machinchose} and \reff{v2point}
\begin{align}\label{machinmachin}
\nonumber\lambda_j\left(v'_j(u)-v'_{j-1}(u-b_j)\right)&=(ru+\lambda_jb_j)v''_j(u)+rv'_j(u)\\
&=(ru+\lambda_jb_j)^{\frac{\lambda_j}{r}-1}\int_u^{\gamma_j}\frac{\lambda_jv''_{j-1}(x-b_j)}{(rx+\lambda_jb_j)^{\frac{\lambda_j}{r}-1}}dx+rv'_j(u).
\end{align}

\vspace{0.5em}
Now if we differentiate \reff{HJB.new} and solve the corresponding ODE for $v'_j$, we obtain
\begin{equation}\label{machintruc}
v'_j(u)=(ru+\lambda_jb_j)^{\frac{\lambda_j}{r}-1}\int_u^{\gamma_j}\frac{\lambda_jv'_{j-1}(x-b_j)}{(rx+\lambda_jb_j)^{\frac{\lambda_j}{r}}}dv-\left(\frac{ru+\lambda_jb_j}{r\gamma_j+\lambda_jb_j}\right)^{\frac{\lambda_j}{r}-1}.
\end{equation}

\vspace{0.5em}
Using \reff{machintruc} in \reff{machinmachin}, we obtain for $u>b_j+b_{j-1}$
\begin{align}\label{condd}
\nonumber&\lambda_j\left(v'_j(u)-v'_{j-1}(u-b_j)\right)\\
\nonumber&=\lambda_j(ru+\lambda_jb_j)^{\frac{\lambda_j}{r}-1}\int_u^{\gamma_j}\frac{ (rx+\lambda_jb_j)v''_{j-1}(x-b_j)+rv'_{j-1}(x-b_j)}{(rx+\lambda_jb_j)^{\frac{\lambda_j}{r}}}dv\\
&\hspace{0.9em}-r\left(\frac{ru+\lambda_jb_j}{r\gamma_j+\lambda_jb_j}\right)^{\frac{\lambda_j}{r}-1}.
\end{align}

\vspace{0.5em}
Then we have for all $x\geq u>b_j+b_{j-1}$ and $x\neq b_j+b_{j-1}+b_{j-2}$
\begin{align*}
 (rx+\lambda_jb_j)v''_{j-1}(x-b_j)+rv'_{j-1}(x-b_j)&=(r(x-b_j)+\lambda_{j-1}b_{j-1})v''_{j-1}(x-b_j)\\
 &\hspace{0.9em}+\left(\lambda_jb_j-\lambda_{j-1}b_{j-1}+rb_j\right)v''_{j-1}(x-b_j)\\
 &\hspace{0.9em}+rv'_{j-1}(x-b_j)\\
 &=\lambda_{j-1}\left(v'_{j-1}(x-b_j)-v'_{j-2}(x-b_j-b_{j-1})\right)\\
 &\hspace{0.9em} +\left(\lambda_jb_j-\lambda_{j-1}b_{j-1}+rb_j\right)v''_{j-1}(x-b_j)\\
 &\leq\left(\lambda_jb_j-\lambda_{j-1}b_{j-1}+rb_j\right)v''_{j-1}(x-b_j),
\end{align*}
where we used the induction hypothesis \reff{hdr} in the last inequality.

\vspace{0.5em}
Since $v_{j-1}$ is concave, the sign of the right-hand side above is given by the sign of
\begin{align*}
\lambda_jb_j-\lambda_{j-1}b_{j-1}+rb_j=\frac{JB}{\qeds}-\frac{(J-1)B}{\qeds}+rb_j=\frac{B}{\qeds}+rb_j\geq 0.
\end{align*}

\vspace{0.5em}
Reporting this in \reff{condd} implies 
$$v'_j(u)-v'_{j-1}(u-b_j)\leq 0,\text{ }u>b_j+b_{j-1}.$$

\vspace{0.5em}
It remains to prove \reff{propz} when $b_j<u<b_j+b_{j-1}$. In that case, \reff{propz} can be written 
$$v'_j(u)-\frac{\overline{v}_{j-1}}{b_{j-1}}\leq 0, \text{ }b_j<u<b_j+b_{j-1},$$
which is equivalent by concavity of $v_j$ to
$$v'_j(b_j^+)-\frac{\overline{v}_{j-1}}{b_{j-1}}\leq 0.$$

\vspace{0.5em}
Now using \reff{machintruc}, we also have
\begin{align}
\nonumber v'_j(b_j^+)&=(ru+\lambda_jb_j)^{\frac{\lambda_j}{r}-1}\int_{b_j}^{b_j+b_{j-1}}\frac{\lambda_j\frac{\overline{v}_{j-1}}{b_{j-1}}}{(rx+\lambda_jb_j)^{\frac{\lambda_j}{r}}}dv\\
\nonumber&\hspace{0.9em}+v'_{j-1}(b_j+b_{j-1})\left(\frac{ru+\lambda_jb_j}{r(b_j+b_{j-1})+\lambda_jb_j}\right)^{\frac{\lambda_j}{r}-1}\\
\nonumber&=\frac{\overline{v}_{j-1}}{b_{j-1}}\frac{\lambda_j}{\lambda_j-r}\left(1-\left(\frac{rb_j+\lambda_jb_j}{r(b_j+b_{j-1})+\lambda_jb_j}\right)^{\frac{\lambda_j}{r}-1}\right)\\
\nonumber&\hspace{0.9em}+v'_{j-1}(b_j+b_{j-1})\left(\frac{ru+\lambda_jb_j}{r(b_j+b_{j-1})+\lambda_jb_j}\right)^{\frac{\lambda_j}{r}-1}.
\end{align}

And thus
\begin{align}
\nonumber v'_j(b_j^+)&\leq\frac{\overline{v}_{j-1}}{b_{j-1}}\frac{\lambda_j}{\lambda_j-r}\left(1-\left(\frac{rb_j+\lambda_jb_j}{r(b_j+b_{j-1})+\lambda_jb_j}\right)^{\frac{\lambda_j}{r}-1}\right)\\
\nonumber&\hspace{0.9em}+v'_{j-1}(b_{j-1}^+)\left(\frac{ru+\lambda_jb_j}{r(b_j+b_{j-1})+\lambda_jb_j}\right)^{\frac{\lambda_j}{r}-1}\\
\nonumber&=\phi_{\frac{b_{j-1}}{b_j}}\left(\frac{r}{\lambda_j}\right)\frac{\overline{v}_{j-1}}{b_{j-1}}\left(\frac{\phi_{\frac{b_{j-1}}{b_j}}\left(\frac{r}{\lambda_j}\right)-1}{\phi_{\frac{b_{j-1}}{b_j}}\left(\frac{r}{\lambda_j}\right)(x-1)}+\frac{v'_{j-1}(b_{j-1}^+)}{\frac{\overline{v}_{j-1}}{b_{j-1}}}\right),
\end{align}
which implies
$$v'_j(b_j^+)-\frac{\overline{v}_{j-1}}{b_{j-1}}\leq\phi_{\frac{b_{j-1}}{b_j}}\left(\frac{r}{\lambda_j}\right)\frac{\overline{v}_{j-1}}{b_{j-1}}\left(\frac{v'_{j-1}(b_{j-1}^+)}{\frac{\overline{v}_{j-1}}{b_{j-1}}}-\psi_{\frac{b_{j-1}}{b_j}}\left(\frac{r}{\lambda_j}\right)\right).$$

By Assumption \ref{decreasing assumption}, we know that $b_j\geq b_{j-1}$, hence with \reff{hyp.lambda} and what we recalled earlier about the functions $\psi_\beta$ in Remark \ref{psi}, we have
$$\frac{\overline{v}_{j-1}}{b_{j-1}}\leq\psi\left(\frac{r}{\lambda_j}\right)\leq\psi_{\frac{b_{j-1}}{b_j}}\left(\frac{r}{\lambda_j}\right),$$
which implies the desired property and ends the proof.

\qed
\end{appendix}


\vspace{0.5em}

\begin{thebibliography}{99}
\bibitem{abreu} Abreu, D., Milgrom, P., Pearce, D. (1991). Information and
timing in repeated partnerships, \textsl{Econometrica}, 59, 1713--1733.

\bibitem{ait} A\"{\i}t-Sahalia, Y., Cacho-Diaz, J., Laeven, R. (2010).
Modeling financial contagion using mutually exciting jump processes. NBER \
Working paper No. 15850.

\bibitem{azizpour} Azizpour, S., Giescke, K. (2008). Self-exciting corporate
defaults: Contagion vs.\ frailty, working paper, Stanford University.

\bibitem{biais2} Biais, B., Mariotti T., Rochet, J.-C., Villeneuve, S.
(2010). Large risks, limited liability and dynamic moral hazard, \textsl{%
Econometrica}, 78(1), 73--118.

\bibitem{bre} Br\'{e}maud, P. (1981). Point processes and queues: martingale
dynamics, Springer Verlag.

\bibitem{cviz} Cvitani\'c, J., Zhang, J. (2010). Contract theory in
continuous time models, monograph, in preparation.

\bibitem{davis} Davis, M., Lo, V. (2001). Infectious defaults, \textsl{%
Quantitative Finance}, 1, 382--387.

\bibitem{dellm} Dellacherie, C., Meyer, P.-A. (1982). Probabilities and
potential, Volume B, Amsterdam, North-Holland.

\bibitem{demarzo4} DeMarzo, P., Fishman, M. (2007a). Agency and optimal
investment dynamics, \textsl{The Review of Financial Studies}, 20, 151--189.

\bibitem{demarzo5} DeMarzo, P., Fishman, M. (2007b). Optimal long-term
financial contracting, \textsl{The Review of Financial Studies}, 20,
2079--2128.

\bibitem{frey} Frey, R., Backhaus, J. (2008). Pricing and hedging of
portfolio credit derivatives with interacting default intensities, \textsl{%
International Journal of Theoretical and Applied Finance}, 11(6), 611-634.

\bibitem{giesecke} Giesecke, K., Kakavand, H., Mousavi, M., Takada, H.
(2010). Exact and efficient simulation of correlated defaults. \textsl{SIAM
Journal on Financial Mathematics}, 1, 868-896.

\bibitem{karat} Karatzas, I., Shreve, S. (1991). Brownian motion and
stochastic calculus, \textsl{Springer-Verlag}, New-York.

\bibitem{jar} Jarrow, R., Yu, F. (2001). Counterparty risk and the pricing
of defaultable securities, \textsl{Journal of Finance}, 53, 2225--2243.

\bibitem{kraft} Kraft, H., Steffensen, M. (2007). Bankruptcy, counterparty
risk, and contagion, \textsl{Review of Finance}, 11, 209--252.

\bibitem{laurent} Laurent, J.-P., Cousin, A., Fermanian, J.-D. (2008).
Hedging default risks of CDOs in Markovian contagion models, preprint.

\bibitem{pages} Pag\`{e}s, H. (2012). Bank monitoring incentives and optimal
ABS, \textsl{Journal of Financial Intermediation}, 10.1016/j.jfi.2012.06.001.

\bibitem{san} Sannikov, Y. (2008). A continuous-time version of the
principal-agent problem, \textsl{Review of Economic Studies}, 75, 957--984.

\bibitem{san2} Sannikov, Y., Skrzypacz, A. (2007). Impossibility of
collusion under imperfect monitoring with flexible production, \textsl{%
American Economic Review}, 97, 1794--1823.

\bibitem{san3} Sannikov, Y., Skrzypacz, A. (2010). The role of information
in repeated games with frequent actions, \textsl{Econometrica}, 78(3),
847--882.

\bibitem{yu} Yu, F. (2007). Correlated defaults in intensity-based models, 
\textsl{Mathematical Finance}, 17(2), 155--173.
\end{thebibliography}
\end{document}